\documentstyle[12pt,amsmath,amssymb]{article}
\begin{document}

\newcommand{\Ext}{{Ext}}
\newcommand{\R}{{\Bbb R}}
\newcommand{\N}{{\Bbb N}}
\newcommand{\Ha}{{\cal H}}
\newcommand{\calS}{{\cal S}}
\newcommand{\Gama}{{\cal F}_{\mathrm \Ext}}
\renewcommand{\kappa}{\varkappa}
\newcommand{\indlim}{\operatornamewithlimits{ind\,lim}}
\newcommand{\supp}{\operatorname{supp}}
\newcommand{\hotimes}{\hat\otimes}
\newcommand{\C}{{\Bbb C}}
\newcommand{\ho}{\hat \otimes}
\newcommand{\la}{\langle}
\newcommand{\ra}{\rangle}
\newcommand{\Fin}{{\cal F}_{{\mathrm fin}}}
\newcommand{\ot}{\otimes}
\newcommand{\id}{\operatorname{id}}
\newcommand{\rom}[1]{{\rm #1}}
\newcommand{\Gamaone}{{\cal F}_{1,k_0}}
\newcommand{\muG}{\mu_{\mathrm G}}
\newcommand{\ltwoG}{(L^2_{\mathrm G})}
\newcommand{\om}{\omega}
\newcommand{\wom}[1]{{:}\,\om^{\ot #1}{:}_{\mathrm G}}
\newcommand{\wup}[1]{{:}\,\upsilon^{\ot #1}{:}_{\mathrm G}}
\newcommand{\Ex}{{\frak X}}
\newcommand{\Pe}{{\cal P}}
\newcommand{\SG}{(S_{\mathrm G})}
\renewcommand{\ll}{\la\!\la}
\newcommand{\rr}{\ra\!\ra}
\newcommand{\wexp}[2]{{:}\,e^{\la #1,#2\ra}\,{:}_{\mathrm G}}
\newcommand{\Hol}{\operatorname{Hol}_0(S_\C)}
\newcommand{\Emin}[1]{{\cal E}^1_{\mathrm min}(#1)}
\newcommand{\chush}{\pmb{\pmb |}}
\newcommand{\Chush}{\pmb{\pmb{\Big|}}}
\newcommand{\ds}{\diamondsuit}
\newcommand{\di}{\partial}
\newcommand{\dig}{\partial^\dag}
\newcommand{\muCP}{\mu_{{\mathrm CP}}}
\newcommand{\vvv}{\phi}
\newcommand{\ttt}{\varphi}
\newcommand{\expp}{e^{\la\om,\ttt\ra}}

\newcommand{\norm}{\pmb{\boldsymbol{\|}}}
\newcommand{\haha}{ \,{|}\!{|}\!{|}\,  }

\renewcommand{\S}{{\cal S}}
\renewcommand{\Box}{\blacksquare}

\renewcommand{\hat}{\widehat}

\newtheorem{lem}{Lemma}
\newtheorem{th}{Theorem}
\newtheorem{prop}{Proposition}

\thispagestyle{empty}

$\mbox{}$\vspace{10mm}
\begin{center}
\LARGE\bf Operators of Gamma White Noise \\[2mm]
\LARGE\bf Calculus\\[20mm]
\large\bf Yuri G. Kondratiev\\[3mm]
\large Inst.\ f.\ Angew.\ Math., Univ.\ Bonn, D-53115 Bonn, Germany; and\\
\large BiBoS, Univ.\ Bielefeld, D-33615 Bielefeld, Germany; and\\
Inst.\ Math., NASU, 252601 Kiev, Ukraine\\[15mm]
\large\bf Eugene W. Lytvynov\\[3mm]
\large Inst.\ f.\ Angew.\ Math., Univ.\ Bonn, D-53115 Bonn, Germany; and\\
BiBoS, Univ.\ Bielefeld, D-33615 Bielefeld, Germany
\vspace{30mm}
\end{center}
\setcounter{section}{-1}
\setcounter{page}{0}

\newpage

\begin{center}\bf Abstract\end{center}
\noindent\begin{small}

The paper is devoted to the study of Gamma white noise analysis.
We define an extended Fock space $\Gama(\Ha)$ over $\Ha=L^2(\R^d,
d\sigma)$, and show how to include the usual Fock space ${\cal F}
(\Ha)$ in it as a subspace. We introduce in $\Gama(\Ha)$
operators $a(\xi)=\int_{\R^d}
dx\,\xi(x)a(x)$, $\xi\in S$, with $a(x)=\dig_x+2\dig_x\di_x+1+\di_x
+\dig_x\di_x\di_x$, where $\dig_x$ and $\di_x$ are the creation and annihilation
operators at $x$. We show that $(a(\xi))_{\xi\in S}$ is a family of
commuting selfadjoint operators in $\Gama(\Ha)$ and construct the Fourier
transform in generalized joint eigenvectors of this family. This transform
is a unitary $I$ between $\Gama(\Ha)$ and the $L^2$-space $L^2(S',d\mu_{\mathrm G})$,
where $\mu_{\mathrm G}$ is the measure of Gamma white noise with intensity $\sigma$.
The image of $a(\xi)$ under $I$ is the operator of multiplication by $\la\cdot,\xi\ra$,
so that $a(\xi)$'s are Gamma field operators. The Fock structure of the Gamma space
determined by $I$ coincides with that discovered in
{\bf [}{\it Infinite Dimensional Analysis,
  Quantum Probability and Related Topics\/} {\bf 1} (1998), 91--117{\bf ]}.
We note that $I$ extends in a natural way the multiple stochastic integral
(chaos) decomposition of the ``chaotic'' subspace of the Gamma space. Next,
we introduce and study spaces of test and generalized functions of Gamma
white noise and derive explicit formulas for the action of the creation, neutral,
and Gamma annihilation operators on these spaces.

\end{small}
\newpage

\tableofcontents
\newpage

\section{Introduction}

During last years, one witnesses the growth of interest to the
the study of compound, or more generally marked Poisson processes
and application of these to different problems of
probability theory, stochastic
analysis and mathematical physics, see e.g.\ \cite{Kallenberg,
Kingman}. The Gamma process is an important example of a compound
Poisson process with a non-finite measure
of the values of jumps of the process in the L\'evy formula.
In this paper, we derive a white  noise calculus for the Gamma
process, studying in detail the so-called Gamma field
operators. Our interest in the Gamma analysis was inspired,
in particular, by the papers by Vershik {\it et al.}\
\cite{Vershik,Vershik1,Vershik2}, where the Gamma measure was used in construction of
a representation of a group of flows. On the other hand, we improve
the results of the paper \cite{Silva}, in which the Fock type
structure of the Gamma space was discovered.

Let us shortly describe our results.
Let $\sigma$
be a non-atomic Radon measure on $\R^d$ and suppose here,
for simplicity, that
$\supp \sigma=\R^d$.
In Sect.~1, we introduce the so-called
extended Fock space $\Gama(\Ha)$ over $\Ha=L^2(\R^d,d\sigma)$.
This space is defined as a
completion of the set ${\cal F}_{\mathrm fin}(S)$ with respect to (w.r.t.)\
a scalar product $(\cdot,\cdot)_{\Gama(\Ha)}$, which is defined by using
a simple combinatorial rule. Here, ${\cal F}_{\mathrm fin}(S):=\bigoplus_{n=0}
^\infty S_\C^{\hat \otimes n}$ is the topological direct sum of the complexified
symmetric tensor powers of the Schwartz test space $S$. Thus,
${\cal F}_{\mathrm fin}(S)$ consists of finite sequences $f=(f^{(0)}, f^{(1)},
\dots,f^{(n)},0,0,\dots)$, where $f^{(j)}$ is a smooth symmetric rapidly
decreasing function of $j$ variables. We show also the way how to include
the usual Fock space over $\Ha$, denoted by ${\cal F}(\Ha)$, into $\Gama(\Ha)$,
so that ${\cal F}(\Ha)$ becomes a proper subspace of $\Gama(\Ha)$.

Recall now that the Gamma white noise measure $\mu_{\mathrm G}$ with intensity
measure $\sigma$ is defined on $S'$---the dual of $S$ w.r.t.\
the zero space $\Ha$---by its Laplace transform
\begin{equation}\label{2.4a}
\ell_{\mathrm G}(\varphi)=
\int_{S'}\exp[\la\om,\varphi\ra]\,d\muG(\om)=\exp\bigg[-\int_{\R^d}\log(1-\varphi(x))\,d\sigma(x)\bigg],
\qquad 1>\varphi\in S.
\end{equation}

In Sect.~2, we introduce a family of commuting, essentially selfadjoint operators
$(a(\xi))_{\xi\in S}$ in $\Gama(\Ha)$ with domain ${\cal F}_{\mathrm fin}(S)$
by the formula
\begin{equation}\label{str}
a(\xi)=a^+(\xi)+2a^0(\xi)+a^-(\xi)+\int_{\R^d}\xi\,d\sigma\cdot\operatorname{id},
\end{equation}
where $a^+(\xi)$ is a usual creation operator:
$$a^+(\xi)\varphi^{\otimes n}=\xi\hat\otimes\varphi^ {\otimes n},$$
$a^0(\xi)$ is a usual neutral operator:
$$a^0(\xi)\varphi^{\otimes n}=n(\xi\varphi)\hat\otimes\varphi^{\otimes(n-1)},$$
and $a^-(\xi)$ is an annihilation operator acting as follows:
$$a^-(\xi)\varphi^{\otimes n}=n\la\xi,\varphi\ra\varphi^{\otimes(n-1)}
+n(n-1)(\xi\varphi^2)\hat\otimes\varphi^{\otimes(n-2)}$$
(notice that the first addend corresponds to the usual annihilation operator,
while the second addend is an annihilation of a new type).

By using the spectral approach to commutative Jacobi fields in the Fock space
\cite{Bee, BeLi,Ly, Benew}, we construct the Fourier transform in generalized joint
eigenvectors of the family $(a^\sim(\xi))_{\xi\in S}$, where $a^\sim(\xi)$
is the closure of $a(\xi)$. This transform, denoted by $I$, is a unitary between
$\Gama(\Ha)$ and the Gamma space $L^2(S',d\mu_{\mathrm G}):=(L^2_{\mathrm G})$,
which was already
constructed in \cite{Silva}. The image of $a^\sim(\xi)$ under $I$ is the operator
of multiplication by the monomial $\la\cdot,\xi\ra$. Thus, $a^\sim(\xi)$'s
are actually  Gamma field operators.

In Sect.~3, we study the ``chaos decomposition'' of the Gamma space which
naturally appears from its Fock type structure. We show also that the restriction
$I\restriction{\cal F}(\Ha)$ of the unitary $I$ to the Fock space as a
subspace of $\Gama(\Ha)$ is exactly the standard isomorphism between
${\cal F}(\Ha)$ and the part of $(L^2_{\mathrm G})$ which
appears as a result of
the chaos
expansion in multiple stochastic integrals w.r.t.\ the compensated Gamma process.
Thus, our approach gives a natural way how, in the Gamma case, to overcome
the main difficulty in compound Poisson analysis connected with the fact that
these processes do not possess the chaotic representation property.

In Sect.~4, we introduce and study spaces of test and generalized functions of Gamma
white noise. This allows us, in particular, to introduce the coordinate operators
$\omega(x)\cdot$, $x\in\R^d$, which, as follows from
\eqref{str}, have the form
$$\omega(x)\cdot=\dig_x+2\dig_x\di_x+1+\di_x+\dig_x\di_x\di_x,$$
where $\dig_x$ and $\di_x$ are the (images of the)  creation and annihilation
operators at $x$.

Finally, in Sects.~5 and 6, we derive explicit formulas for the action of the
operators $a^+(\xi)$, $a^0(\xi)$, and $a^-(\xi)$. Some formulas appearing here
are nothing but infinite-dimensional analogs of the formulas obtained in the
one-dimensional case by Meixner in his classical work \cite{Meixner}.

\section{Extended Fock space and its rigging}

Let $\sigma$ be a Borel, regular, non-atomic,  $\sigma$-finite
measure on
$\R^d$, $d\in \N$. First, we
recall the construction from \cite{Hir} of a rigging
 of the real
$L^2$-space
$L^2(\R^d,d\sigma)=\Ha$ by  spaces of test and generalized functions of
Schwartz type. (Notice that, since $\sigma$ is not necessarily the
Lebesgue measure,
we need some additional consideration.)

Let $(e_j)_{j=0}^\infty$ be the system of Hermite functions on $\R$. For each
$p\ge1$, define ${\cal S}_p(\R)$ to be the real Hilbert space spanned by the
orthonormal basis $(e_j(2j+2)^{-p})_{j=0}^\infty$, and let ${\cal S}_p
(\R^d)={\cal S}_p(\R)^{\otimes d}$.
Considered as a subspace of $L^2(\R^d,dx)$, every space ${\cal S}_p(\R^d)$
coincides with the domain of the operator $(H^{\ot d})^p$, where $H^{\ot d}$ is
the harmonic oscillator in $L^2(\R^d,dx)$: $H^{\ot d}=-\sum_{i=1}^d\big(\frac
d{dx_i}\big)^2+\sum_{i=1}^d x_i^2+1$.
As well known, $\calS(\R^d)=\projlim_{
p\to\infty}\calS_p(\R^d)$ is the Schwartz space of rapidly decreasing functions
on $\R^d$, $\calS_1(\R^d)$
consists of continuous functions, and
$$ \R^d\ni x\mapsto \delta_x\in \calS_{-1}(\R^d)$$
is a continuous mapping, where $\calS_{-p}(\R^d)$ denotes the dual of
$\calS_p(\R^d)$.
We suppose the existence of $\epsilon\ge0$ such that
the space $\calS_{1+\epsilon}(\R^d)$ is continuously embedded into
$\Ha$; for example, the following condition holds:
$$\int_{\R^d}\|\delta_x\|^2_{\calS_{-1-\epsilon}(\R^d)}\,d\sigma(x)
<\infty.$$
Let $O_p\colon S_p(\R^d)\hookrightarrow\Ha$ be an embedding operator.
Since $\calS(\R^d)$ is a nuclear space, we can suppose without loss of
generality that the operator
$O_{1+\epsilon}$  is of Hilbert--Schmidt type
(in case of the Lebesgue measure, $d\sigma(x)=dx$, we can take
$\epsilon=0$).
Note that, because of the regularity and $\sigma$-finiteness of the
measure $\sigma$, $\calS(\R^d)$ is a dense subset of $L^2(\R^d,\sigma)$.

Define now $S_p$ to be the Hilbert factor space $\calS_{p+\epsilon}/\ker
O_{p+\epsilon}$. By \cite{BeKo}, Ch.~5, Sect.~5, subsec.~1,
$\big\{\Ha,\, S_p\mid p\ge1\big\}$ is a sequence of compatible Hilbert spaces.
Thus, we
obtain the rigging
\begin{equation}\label{1.1}
S'=\indlim_{p\to\infty} S_{-p}\supset
L^2(\R^d,d\sigma)=\Ha\supset\projlim_{p\to\infty} S_p=S,
\end{equation}
where $S_{-p}$, resp.\ $S'$ is the dual of $S_p$, resp.\ $S$
w.r.t.\ the zero space $\Ha$.
We stress that, for arbitrary $p\ge p'$, the space $S_p$ is continuously
embedded into $S_{p'}$ and $|\cdot|_{p}\ge|\cdot|_{p'}$, where $|\cdot|_p$
denotes the $S_p$ norm.
Notice that, in fact,  the spaces $S_p$ and $S$ are completely
determined by the support of $\sigma$, i.e., $S$
is the Schwartz test space on $\supp\sigma$.

Now we wish to define an $n$-particle extended Fock space over $\Ha$,
$\Gama^{(n)}(\Ha)$, $n\in \N$. Under a loop $\kappa$
connecting points $x_1,\dots,x_m$, $m\ge2$, we understand a class of ordered
sets $(x_{\pi(1)},\dots,x_{\pi(m)})$, where $\pi$ is a permutation of
$\{1,\dots,n\}$, which coincide up to a cyclic permutation. For example,
$(x_1,x_2,x_3)$ and $(x_3,x_1,x_2)$ define the same loop, while
$(x_1,x_2,x_3)$ and $(x_2,x_1,x_3)$ define different loops. We put
$|\kappa|=m$.
We will interpret also a set $\{x\}$ as a ``one-point'' loop $\kappa$, i.e., a loop that
comes out of $x$, $|\kappa|=1$.

Let $\alpha_n=\{\kappa_1,\dots,\kappa_{|\alpha_n|}\}$ be a collection of loops
$\kappa_j$ that connect points from the set $\{x_1,\dots,x_n\}$ so that every
point $x_i\in\{x_1,\dots,x_n\}$ goes into one loop $\kappa_j=\kappa_{j(i)}$
from $\alpha_n$. Here, $|\alpha_n|$ denotes the number of the loops in
$\alpha_n$, evidently $n=\sum_{j=1}^{|\alpha_n|}|\kappa_j|$.

Let $A_n$ stand for the set of all possible collections of loops $\alpha_n$
over the points $\{x_1,\dots,x_n\}$. Every $\alpha_n\in A_n$ generates the
following continuous mapping
\begin{multline}\label{1.2}
 S^{\hotimes n}_\C\ni f^{(n)}=f^{(n)}(x_1,\dots,x_n)\mapsto\\
\mapsto f^{(n)}_{\alpha_n}(
\underbrace{x_1,\dots,x_1}_{\text{$|\kappa_1|$ times}} ,
\underbrace{x_2,\dots,x_2}_{\text{$|\kappa_2|$ times}} ,\dots,
\underbrace{x_{|\alpha_n|},
\dots,x_{|\alpha_n|}}_{\text{$|\kappa_{|\alpha_n|}|$ times}} )\in
S_\C^{\otimes|\alpha_n|},
           \end{multline}
where the lower index $\C$ denotes complexification of a real space and
the symbol $\hotimes$ stands for the symmetric
tensor power. Indeed, for any $p\ge1$, the diagonalization operator $\frak D
$ given by
$$f^{(2)}=f^{(2)}(x_1,x_2)\mapsto {\frak D}f^{(2)}=\big({\frak
  D}f^{(2)}\big)(x)=f^{(2)}(x,x)$$
acts continuously from $S_{p,\C}^{\hotimes 2}$ into $S_{p,\C}$,
\begin{equation}\label{1.3}
|{\frak D}f^{(2)}|_p\le C_p|f^{(2)}|_p,\end{equation} where
$|\cdot|_p$ denotes also the norm of each space $S_{p,\C}^{\otimes
n}$, $n\in \N$, which yields
\begin{equation}\label{1.4}
|f^{(n)}_{\alpha_n}|_p\le C_p^{(|\kappa_1|-1)+(|\kappa_2|-1)+\cdots+
(|\kappa_{|\alpha_n|}|-1)}|f^{(n)}|_p,\end{equation}
giving the continuity of the mapping \eqref{1.2}.

Thus, we define a scalar product on $S_{\C}^{\hotimes n}$ by
\begin{equation}\label{1.5}
( f^{(n)},g^{(n)})_{\Gama^{(n)}(\Ha)}=\sum_{\alpha_n\in A_n}\int
_{\R^{d|\alpha_n|}}\big(\overline{ f^{(n)}}g^{(n)}\big)_{\alpha_n}\,
d\sigma^{|\alpha_n|},
\end{equation}
where $\overline{f^{(n)}}$ is the complex conjugate of $f^{(n)}$. Let
$\Gama^{(n)}(\Ha)$ be
the closure of $S_{\C}^{\ho n}$ in the norm generated by \eqref{1.5}.
\vspace{2mm}

\noindent {\it Remark\/} 1. It is easy to see that the number of the summands in
the series \eqref{1.5} is exactly $n!$. Let us prove this by induction. For
$n=1$, this is trivial. Let the statement hold for $n$. At the step
$(n+1)$ we
add a point $x_{n+1}$ to the set $
\{x_1,\dots,x_n\}$. If $x_{n+1}$ goes into a ``one-point'' loop, then by the
supposition there are exactly $n!$ variants. If $x_{n+1}$ goes into a loop that
connect this point with other points from $\{x_1,\dots,x_n\}$, then to each
loop connecting $i$
points from $\{x_1,\dots,x_n\}$ we can add the point $x_{n+1}$ in $i$ different
ways. That gives $n\,n!$ variants. Thus, in total we have exactly $(n+1)!$
variants.\vspace{2mm}

The extended Fock space $
\Gama(\Ha)$
over $\Ha$ is defined as a weighted direct sum of the spaces $\Gama^{(n)}(\Ha)$:
\begin{equation}\label{1.5a}
\Gama(\Ha)=\bigoplus_{n=0}^\infty \Gama^{(n)}(\Ha)\,n!,\end{equation}
where $\Gama^{(0)}(\Ha)=\C$ and  $0!=1$. I.e., $\Gama(\Ha)$ consists of
sequences $f=(f^{(0)}, f^{(1)},$  $f^{(2)},\dots)$
such that $f^{(n)}\in \Gama^{(n)}(\Ha)$ and
$$ \|f\|_{\Gama(\Ha)}^2=\sum_{n=0}^\infty\|f^{(n)}\|^2_{\Gama^{(n)}(\Ha)}
n!<\infty.$$
Throughout the paper we will identify $f^{(n)}\in \Gama^{(n)}(\Ha)$ with
the vector
$$(0,\dots,0,f^{(n)},0,0\dots)\in\Gama(\Ha).$$

In the following section, we will need
a nuclear space that is topologically, i.e., densely and continuously, embedded
into $\Gama(\Ha)$, and its dual space. Thus, we take
\begin{equation}
\label{a}
\Gama(\Ha)\supset \Fin(S),\qquad \Fin^*(S)\supset\Fin(S).
\end{equation}
Here, $\Fin(S)$ is the topological direct sum of the spaces
$S_{\C}^{\ho n}$, i.e., $\Fin(S)$ consists of all finite sequences
$f=(f^{(0)},f^{(1)},\dots,f^{(m)},0,0,\dots)$ such that $f^{(n)}\in S_\C^{\ho
  n}$ and the convergence in $\Fin(S)$ means the uniform finiteness and the
coordinate-wise convergence. That $\Fin(S)$ is dense in $\Gama(\Ha)$
follows directly from the definition of $\Gama(\Ha)$. \eqref{1.4}, \eqref{1.5},
and
Remark~1 give
\begin{align}
\|f^{(n)}\|^2_{\Gama^{(n)}(\Ha)}&\le n!\max_{\alpha_n\in A_n}\|f^{(n)}_{\alpha_n}
\|^2_{\Ha_\C^{\ho |\alpha_n|}}\notag\\
&\le n!K_1^{2n}\max_{\alpha_n\in A_n}|f^{(n)}_{|\alpha_n|}|_1^2\le
n!\,K_1^{2n}C_1^{2(n-1)}|f^{(n)}|_1^2, \label{1.6}\end{align}
where $K_1$ is the norm of the inclusion operator $S_1\hookrightarrow\Ha$ (we suppose that
the constants $K_1,\,C_1\ge1$). \eqref{1.6} implies the continuity of the
embedding $\Fin(S)\hookrightarrow\Gama(\Ha)$.

The space $\Fin^*(S)$ in \eqref{a} is the dual of $\Fin(S)$.
It will be convenient for us to take it as the dual of $\Fin(S)$ w.r.t.\ the
zero space ${\cal F}(\Ha)$, the usual Fock space over $\Ha$:
$${\cal F}(\Ha)=\bigoplus_{n=0}^\infty \Ha_\C^{\ho n}\,n!.$$
So, the second inclusion in \eqref{a} is  part of the nuclear triple (cf.\ e.g.\
\cite{BeKo})
\begin{equation}\label{standard}
\Fin^*(S)\supset{\cal F}(\Ha)\supset\Fin(S).\end{equation}
The space $\Fin^*(S)$ consists of infinite sequences
$F=(F^{(0)},F^{(1)},F^{(2)},\dots)$, where $F^{(n)}
\in S_\C^{\prime\,\ho n}$ and the dualization with $f\in\Fin(S)$ is given by
$$ \la\!\la F,f\ra\!\ra=\sum_{n=0}^\infty \la \overline{F^{(n)}},f^{(n)}
\ra\,n!,$$
where $\la\cdot,\cdot\ra$ denotes the dualization generated
by the scalar product in $\Ha^{\ho n}$, which is supposed to be
linear in  both dots.

Finally, we note that the usual Fock space ${\cal F}(\Ha)$
can be included into $\Gama(\Ha)$. To this end, we construct a subset of $S_\C^{\ho
  n}$
of the form
$${\cal R}^{(n)}=\operatorname{l.s.}\big\{\,
\ttt_1\ho\cdots\ho\ttt_n\mid\ttt_i\in
S_\C,\
\sigma (\supp \ttt_i\cap\,\supp\ttt_j)=0\ \forall i,j=1,\dots,n,\, i\ne j\},$$
where $\operatorname{l.s.}$ denotes the linear span.
It follows from \eqref{1.5} that, for arbitrary $f^{(n)},g^{(n)}\in{\cal
  R}^{(n)}$,
\begin{equation}\label{1.7}
( f^{(n)},g^{(n)})_{\Gama^{(n)}(\Ha)}=( f^{(n)},g^{(n)})_{\Ha_\C^{\ho n}}.\end{equation}
Let $\overline{{\cal R}^{(n)}}$
be the closure of ${\cal R}^{(n)}$ in the $\Gama^{(n)}(\Ha)$ norm. Evidently,
$\overline{{\cal R}^{(n)}}$ is a subspace of $\Gama^{(n)}(\Ha)$. On the other
hand, the set ${\cal R}^{(n)}$ is dense in $\Ha_\C^{\ho n}$. Therefore, in view
of \eqref{1.7}, we can identify $\overline{{\cal R}^{(n)}}$ with $\Ha_\C^{\ho
  n}$.
Then, the subspace $\bigoplus_{n=0}^\infty \overline{{\cal R}^{(n)}}\,n!$ of $\Gama(\Ha)$
is identified with ${\cal F}(\Ha)$.\vspace{2mm}

\noindent {\it Remark\/} 2. Let us explain how the extended Fock space happens to
be greater than the usual Fock space. Take, for example, $n=2$. Then, the space
$\Ha_\C^{\ho2}=\hat L^2(\R^{2d},d\sigma^2)$ is the
complex space of quadratic integrable symmetric functions of two variables from $\R^d$. Of
course, any function $f^{(2)}$ from $\Ha_\C^{\ho 2}$ is defined on $\R^{2d}$ up
to a set of zero $\sigma^2$ measure. The diagonal $\{(x,x)\mid x\in \R^d\}$ in
$\R^{2d}$ has $\sigma^2$ measure zero, while $\sigma$ is non-atomic. On the
other hand, every function $f^{(2)}\in\Gama^{(2)}(\Ha)$ must also be defined on
the diagonal up to a set of zero $\sigma$ measure, if we identify this diagonal
with $\R^d$. Thus, any function $f^{(2)}\in\Ha_\C^{\ho 2}$ is included into
$\Gama^{(2)}(\Ha)$ if we put additionally that $f^{(2)}$ is equal to zero on the
diagonal. An analogous situation takes place in all the other $n$-particle
spaces with $n\ge3$.\vspace{2mm}

\noindent{\it Remark\/} 3. Let us stress that we have constructed  two
different types of inclusion: the first one $\Fin(S)\subset{\cal F}(\Ha)$ in
\eqref{standard} and ${\cal F}(\Ha)$ as a subspace of $\Gama(\Ha)$.

\section{Gamma field operators\\ and the Fourier transform}

In this section, we will define a family of Gamma field operators and construct
a Fourier transform in generalized joint eigenvectors of this family.

For each $\xi\in S$, let $a^+(\xi)$ be the standard creation operator defined
on $\Fin(S)$:
\begin{equation}\label{2.0}
a^+(\xi)f^{(n)}=\xi\ho f^{(n)},\qquad f^{(n)}\in S^{\ho n}_\C,\ n\in
\N_0=\{0,1,2,\dots\}.
\end{equation}
Because of the estimate
\begin{equation}\label{2.1}
|\xi\ho f^{(n)}|_p\le|\xi|_p|f^{(n)}|_p,\end{equation}
the operator $a^+(\xi)$ acts continuously on $\Fin(S)$.

Let us calculate the adjoint operator of $a^+(\xi)$ in $\Gama(\Ha)$, which will
be denoted by $a^-(\xi)$. By \eqref{1.5}, \eqref{1.5a}, and \eqref{2.0}, we have for arbitrary
$\xi,\ttt,\psi\in S$,
\begin{gather*}
( a^+(\xi)\ttt^{\otimes
  n},\psi^{\otimes(n+1)})_{\Gama(\Ha)}=(n+1)!\,(\xi\ho\ttt^{\otimes
  n},\psi^{\otimes(n+1)}
)_{\Gama^{(n+1)}(\Ha)}\\
=(n+1)!\,\sum_{\alpha_{n+1}\in
  A_{n+1}}\int_{\R^{d|\alpha_{n+1}|}}\big((\xi\psi)\ho(\ttt\psi)^{\otimes
  n}\big)_{\alpha_{n+1}}
\,d\sigma^{|\alpha_{n+1}|}
\\
=(n+1)!\,\bigg[\la\xi,\psi\ra\sum_{\alpha_n\in
  A_n}\int_{\R^{d|\alpha_n|}}\big((\ttt\psi)^{\otimes
  n}\big)_{\alpha_n}\,d\sigma^{|\alpha_n|}
\\
\text{}
+n\sum_{\alpha_n\in
  A_n}\int_{\R^{d|\alpha_n|}}\big((\xi\ttt\psi^2)\ho(\ttt\psi)^{\otimes(n-1)}\big)
_{\alpha_n}\,d\sigma^{|\alpha_n|}\bigg]\\
=\big(\ttt^{\otimes n},(n+1)\la\xi,\psi\ra\psi^{\otimes
  n}+(n+1)n(\xi\psi^2)\ho\psi
^{\otimes(n-1)}\big)_{\Gama(\Ha)},
\end{gather*}
where we used just as in Remark~1 the observation that a new point can be added
to a loop connecting $i$ points in $i$ different ways. Thus,
$$ a^-(\xi)=a_1^-(\xi)+a_2^-(\xi),$$
where $a_1^-(\xi)$ is the standard annihilation operator:
$$a_1^-(\xi)\ttt^{\otimes n}=n\la\xi,\ttt\ra\ttt^{\otimes (n-1)}$$
and $a_2^-(\xi)$ given by
$$a_2^-(\xi)\ttt^{\otimes n}=n(n-1)(\xi\ttt^2)\ho\ttt^{\otimes(n-2)}$$
is a ``Gamma annihilation'' operator, which appears because of the nonstandard
scalar product in our Fock space. Due to \eqref{1.3}, $a^-(\xi)$ acts
continuously on $\Fin(S)$ and
\begin{align}
|a_1^-(\xi)f^{(n)}|_p&\le n|\xi|_{-p}|f^{(n)}|_p,\notag\\
|a_2^-(\xi)f^{(n)}|_p&\le n(n-1)C_p^2|\xi|_p|f^{(n)}|_p,\label{2.2}
\end{align}

Finally, we define on $\Fin(S)$ the neutral operator $a^0(\xi)$, $\xi\in S$,
as
the differential second quantization of the operator of multiplication by
$\xi$:
$$a^0(\xi)\ttt^{\otimes n}=n(\xi\ttt)\ho\ttt^{\otimes(n-1)},\qquad \ttt\in S.$$
Again, $a^0(\xi)$ acts continuously on $\Fin(S)$ and
\begin{equation}\label{2.3}
|a^0(\xi)f^{(n)}|_p\le nC_p|\xi|_p|f^{(n)}|_p,\qquad f^{(n)}\in S_\C^{\ho
  n}.\end{equation}

Thus, we are in position to define, for each $\xi\in S$, the Gamma field
operator $a(\xi)$ on $\Fin(S)$:
\begin{equation}\label{2.4}
a(\xi)=a^+(\xi)+2a^0(\xi)+\la\xi\ra
\id+a^-(\xi),\end{equation}
where $\la\xi\ra=\int_{\R^d}\xi(x)\,d\sigma(x)$ and $\id$ denotes the identity
operator. Each $a(\xi)$ with domain $\Fin(S)$ is a Hermitian operator in
$\Gama(\Ha)$.

\begin{lem}
The operators $a(\xi)$\rom, $\xi\in S$\rom, with domain $\Fin(S)$ are
essentially selfadjoint in $\Gama(\Ha)$ and their closures
$a^\sim(\xi)$ constitute a family of commuting selfadjoint operators\rom, where
the commutation is understood in the sense of the resolutions of the
identity\rom.
\end{lem}

\noindent {\it Proof}. Let us show
that every $f^{(n)}\in S_\C^{\ho n}$, $n\in \N_0$, is an analytical vector of
each $a(\xi)$, i.e., the series
$$\sum_{m=0}^\infty \frac{\|a(\xi)^mf^{(n)}\|_{\Gama(\Ha)}}{m!}\,|z|^m,\qquad
z\in\C,$$
has a positive radius of convergence.

To this end, we define, for each $p\ge1$ and $k\in\N_0$, a Hilbert space
\begin{equation}\label{space}{\cal F}_{1,k}(S_p)=\bigoplus_{n=0}^\infty S_{p,\C}^{\ho n}\,(n!)^2
2^{nk},\end{equation}
i.e., for $f=(f^{(0)},f^{(1)},f^{(2)},\dots)\in {\cal F}_{1,k}(S_p)$,
$$\|f\|^2_{{\cal F}_{1,k}(S_p)}=\sum_{n=0}^\infty
|f^{(n)}|^2_p\,(n!)^2 2^{nk}.$$
Because of \eqref{1.6}, we have that topologically
${\cal F}_{1,k_0}(S_1)\subset \Gama(\Ha)$ if $2^{k_0/2}\ge K_1C_1$
and
$$\|\cdot\|_{{\cal F}_{1,k_0}(S_1)}\ge\|\cdot\|_{\Gama(\Ha)}.$$
Taking to notice that each $a(\xi)^mf^{(n)}$ belongs to $\Fin(S)$, it suffices
to prove that the series
$$\sum_{m=0}^\infty\frac{\|a(\xi)^mf^{(n)}\|_{{\cal F}_{1,k_0}(S_1)}}{m!}\,|z|^m$$
has a positive radius of convergence. We have
\begin{align*}
a(\xi)^m&=(a^+(\xi)+2a^0(\xi)+\la\xi\ra\id+a^-_1(\xi)+a_2^-(\xi))^m\\
&=\sum_{j=1}^{5^m}Y_{j,m}(\xi)\end{align*}
where every operator $Y_{j,m}(\xi)$ is a composition of $m$ operators every one
of which is one of the operators
$a^+(\xi)$, $2a^0(\xi)$, $\la\xi\ra\id$, $a^-_1(\xi)$, $a_2^-(\xi)$. From
\eqref{2.1}--\eqref{2.3}, we deduce that
\begin{align*}
\|a^+(\xi)f^{(n)}\|_{\Gamaone(S_1)}&\le(n+1)2^{k_0/2}|\xi|_1\|f^{(n)}\|_{\Gamaone(S_1)},\\
\|2a^0(\xi)f^{(n)}\|_{\Gamaone(S_1)}&\le2nC_1|\xi|_1\|f^{(n)}\|_{\Gamaone(S_1)},\\
\|a_1^-(\xi)f^{(n)}\|_{\Gamaone(S_1)}&\le 2^{-k_0/2}|\xi|_{-1}\|f^{(n)}\|_{\Gamaone(S_1)},\\
\|a_2^-(\xi)f^{(n)}\|_{\Gamaone(S_1)}&\le
(n-1)C_1^2 2^{-k_0/2}|\xi|_1\|f^{(n)}\|_{\Gamaone(S_1)}
.\end{align*}
Hence,
$$\|a(\xi)^mf^{(n)}\|_{\Gamaone(S_1)}\le
5^m(2C_1^2 2^{k_0/2})^m\frac{(n+m)!}{n!}\,\max\{|\xi|_1,\la\xi\ra\}
\|f^{(n)}\|_{\Gamaone(S_1)}.$$
Thus, it remains to note that the series
$$\sum_{m=0}^\infty\frac{(10 C_1^2 2^{k_0/2})^m(n+m)!}{m!}\,|z|^m $$
converges as $|z|<(10 C_1^2 2^{k_0/2})^{-1}$.
By using, e.g., \cite{BeKo}, Ch.~5, Th.~1.7, we conclude that the operators
$a(\xi) $ are essentially selfadjoint on $\Fin(S)$.

Let us show that the operators $a(\xi)$ commute on $\Fin(S)$.
Any operators $a^+(\xi_1)$ and $a^+(\xi_2)$ evidently
commute. Hence, their adjoints $a^-(\xi_1)$ and $a^-(\xi_2)$ also commute. Next,
the operators of second quantization $a^0(\xi_1)$ and $a^0(\xi_2)$ commute,
 since
do the operators of multiplication by $\xi_1$ and $\xi_2$.  Next, we have to
show that
\begin{equation}\label{star}
a^+(\xi_1)a^0(\xi_2)+a^0(\xi_1)a^+(\xi_2)=a^+(\xi_2)a^0(\xi_1)+a^0(\xi_2)a^+(\xi_1),
\end{equation}
which evidently yields
$$
a^-(\xi_1)a^0(\xi_2)+a^0(\xi_1)a^-(\xi_2)=a^-(\xi_2)a^0(\xi_1)+a^0(\xi_2)a^-(\xi_1).$$
But \eqref{star} can be easily verified on the vectors of the form $\ttt^{\ot
  n}$, $\ttt\in S$, if one uses the equality
\begin{align*}
a^0(\xi)\ttt_1\ho\ttt_2\ho\cdots\ho\ttt_n&=(\xi\ttt_1)\ho\ttt_2\ho\cdots\ho\ttt_n\\
&\quad+\ttt_1\ho(\xi\ttt_2)\ho\cdots\ho\ttt_n+\cdots+\ttt_1\ho\cdots\ho\ttt_{n-1}\ho(\xi\ttt_n),
\end{align*}
which gives, in particular,
\begin{equation}\label{starstar}
a^0(\xi)\ttt_1\ho\ttt_2^{\ot(n-1)}
=(\xi\ttt_1)\ho\ttt_2^{\ot(n-1)}+(n-1)\ttt_1\ho(\xi\ttt_2)\ho\ttt_2^{\ot(n-2)}.
\end{equation}
Finally, we have to show that
$$a^+(\xi_1)a^-(\xi_2)+a^-(\xi_1)a^+(\xi_2)=a^+(\xi_2)a^-(\xi_1)+a^-(\xi_2)a^+(\xi_1).$$
It is not hard to see that
$$ a_2^-(\xi)\ttt_1\ho\ttt_2^{\ot n}=2n(\xi\ttt_1\ttt_2)\ho\ttt_2^{\ot(n-1)}+n
(n-1)\ttt_1\ho(\xi\varphi_2^2)\ho\varphi_2^{\ot (n-2)}.$$
Then,
\begin{align*}
&(a^+(\xi_1)a^-(\xi_2)+a^-(\xi_1)a^+(\xi_2))\varphi^{\ot n}\\
&\qquad
=n\la\xi_2,\varphi\ra\xi_1\ho\varphi^{\ot(n-1)}+n(n-1)\xi_1\ho(\xi_2\varphi^2)\ho\varphi^{\ot(n-2)}
+\la\xi_1,\xi_2\ra\varphi^{\ot
  n}\\
&\qquad\quad+n\la\xi_1,\varphi\ra\xi_2\ho\varphi^{\ot(n-1)}
+
2n(\xi_1\xi_2\varphi)\ho\varphi^{\ot(n-1)}+n(n-1)\xi_2\ho(\xi_1\varphi^2)\ho\varphi^{\ot
  (n-2)},
\end{align*}
which is symmetric in $\xi_1$ and $\xi_2$.
Thus, arbitrary $a(\xi_1)$ and $a(\xi_2)$ commute.

Since the operators $a^\sim(\xi)$ are essentially selfadjoint on $\Fin(S)$, the
set
$$(a^\sim(\xi)-z\id)\Fin(S),\qquad z\in\C,\ \Im z\ne 0,$$
is dense in $\Gama(\Ha)$. In addition
$$(a^\sim(\xi)-z\id)\Fin(S)\subset\Fin(S).$$
Therefore, for arbitrary $\xi_1,\xi_2\in S$, the operators $a^\sim(\xi_1)$,
$a^\sim(\xi_2)$, and
$$a^\sim(\xi_1)\restriction(a^\sim(\xi_2)-z\id)\Fin(S)$$
have a total
set of analytical vectors. Thus, by \cite{BeKo}, Ch.~5, Th.~1.15, the operators
commute in the sense of the resolutions of the identity. \quad $\Box$

\begin{th} There exists a unitary isomorphism
$$ I\colon \Gama(\Ha)\to L^2(S',{\cal B}(S'), d\muG)=\ltwoG$$
between the extended Fock space $\Gama(\Ha)$ and the complex $L^2$-space over $S'$
with the Borel measure $\muG$ of Gamma white noise\rom, whose Laplace transform
is given by \eqref{2.4a}\rom.
This isomorphism is defined on the dense set $\Fin(S)$ by the formula
\begin{equation}\label{2.5}
\Fin(S)\ni f=(f^{(n)})_{n=0}^\infty\mapsto If=(If)(\om)=\sum_{n=0}^\infty
\la\wom n,f^{(n)}\ra\end{equation}
\rom(the series in \eqref{2.5} is actually finite\rom)\rom, where $\wom n\in
S^{'\,\ho n}$ is the $n$-th Gamma-Wick power of $\om\in S'$ given by the recurrence
formula
\begin{equation}\label{2.6}\begin{gathered}
\wom{(n+1)}=\wom{(n+1)}(x_1,\dots,x_{n+1})=\big(\wom{n}(x_1,\dots,x_n)\om(x_{n+1})\big)^\sim\\
\text{}-n\big(
\wom {(n-1)}(x_1,\dots,x_{n-1})\delta(x_{n+1}-x_n)
\big)^\sim\\
\text{}-n(n-1)\big(
\wom{(n-1)}(x_1,\dots,x_{n-1})\delta(x_{n}-x_{n-1})\delta(x_{n+1}-x_n)
\big)^\sim\\
\text{}-2n\big(
\wom n(x_1,\dots,x_n)\delta(x_{n+1}-x_n)
\big)^\sim
-\big(
\wom n(x_1,\dots,x_n)1(x_{n+1})
\big)^\sim,
\\
\wom0=1,\ \wom1=\om-1,
\end{gathered}\end{equation}
where $(\cdot)^\sim$ denotes the symmetrization of a function.

The image of each
operator $a^\sim(\xi)$\rom, $\xi\in S$\rom, under $I$ is the operator of
multiplication by the monomial $\la\cdot,\xi\ra$ in $\ltwoG$\rom:
\begin{equation}\label{2.7}
Ia^\sim(\xi)I^{-1}=\la\cdot,\xi\ra\cdot\, .
\end{equation}
\end{th}

\noindent {\it Remark\/} 4. Let us explain the term ``Gamma noise,''
introduced in \cite{Silva}.
Let $\Delta$
be a bounded Borel set in $\R^d$ and let $\chi_\Delta$ denote the indicator of
$\Delta$. Then, one can define the random variable (r.v.)\
$X_\Delta=X_\Delta(\om)=\la\om, \chi_\Delta\ra$ as an $\ltwoG$-limit of a
sequence of r.v.'s $\la\om,\xi_j\ra$ such that $\xi_j\to\chi_\Delta$ in
$L^2(\R^d,d\sigma)$ as $j\to\infty$. Then, the Laplace transform
$$\int_{S'}\exp[\lambda X_\Delta(\om)]\,d\muG(\om)=\exp\big[
-\sigma(\Delta)\log(1-\lambda)
\big],\qquad \lambda<1,$$
coincides with the Laplace transform of a r.v.\ having Gamma distribution
with parameter $\sigma(\Delta)$, i.e.,
the density of the distribution of $X_\Delta$ is
\begin{equation}\label{gamma}
p_\Delta(s)=\frac{s^{\sigma(\Delta)-1}e^{-s}}{\Gamma(\sigma(\Delta))},\qquad
s>0.\end{equation}
In case $d=1$ and $d\sigma(x)=dx$, the process $(X_t)_{t\ge0}$, where
$X_t=X_{[0,t]}$ and $X_0=0$, is known as a Gamma process, see, e.g., \cite{Ta}, Sect.~19.
Thus, the triple $(S',{\cal B}(S'),\muG)$ is a direct representation of the
generalized stochastic process $(\dot X_t)_{t\ge0}$, see \cite{GeVi}. So, the term ``Gamma
noise'' is natural for $\muG$.

Let us also note that the Gamma white noise is a special case of
a compound
Poisson white noise, see e.g.\  \cite{GGV75,GeVi,LyReShch,Silva},
because
the Laplace transform \eqref{2.4a} can
also be presented in the form
\begin{equation}\label{CP}
\int_{S'}\exp[\la\om,\varphi\ra]\,d\muG(\om)
=\exp\bigg[\int_0^\infty\!\!\int_{\R^d}(e^{s\varphi(x)}-1)\,d\sigma(x)\,\frac{e^{-s}}{s}\,
ds\bigg].
\end{equation}

This yields, in particular, that the measure $\muG$ is concentrated on the
following subset of $S'$, e.g., \cite{Kal83,KMM78}. Let $\Gamma$
denote the marked configuration space over
$\supp\sigma$ with marks in $]0,\infty[$ that is defined as follows:
$\Gamma$ consists of sets $\gamma\subset \supp\sigma\times]0,\infty[$
such that, for arbitrary disjoint $(x,s_x),(y,s_y)\in\gamma$, we have $x\ne y$, and
for any $a>0$ and $K\subset\R^d$ compact, $\gamma\cap(K\times[a,\infty[)$ is a
finite set.
 Then, the measure $\muG$ is
concentrated on the set of those $\omega\in S'$ which can be presented in the
form
\begin{equation}\label{configuration}
\omega=\sum_{(x,s_x)\in\gamma}s_x\delta_x\end{equation}
with $\gamma\in\Gamma$.\vspace{2mm}

\noindent {\it Proof of Theorem\/} 1. Just as in the case of Gaussian and
Poisson measures, see \cite{Bee,BeLi,Ly}, we deduce, by using the projection spectral theorem
\cite{BeKo}, Ch.~3, Th.~2.7, the existence of a unitary isomorphism $I$ between
$\Gama(\Ha)$ and an $L^2$-space $L^2(S',{\cal B}(S'),d\mu)$ for which the
formulas \eqref{2.5}--\eqref{2.7} hold. Here, $\mu$ is a Borel probability
measure on $S'$, which is the spectral measure of the family
$(a^\sim(\xi))_{\xi\in S}$. Indeed, the major step in the proof of this fact is the
following lemma.

\begin{lem}
\rom{1)} The family $(a^\sim(\xi))_{\xi\in S}$ is connected with
\eqref{a} in a standard way\rom, i\rom.e\rom{.,} for each $\xi\in S$\rom, $a(\xi)$
is a linear continuous operator on $\Fin(S)$\rom.

\rom{2)} For an arbitrary fixed $f\in\Fin(S)$\rom, the mapping
$$ S\ni\xi\mapsto a(\xi)f\in\Fin(S)$$
is linear and continuous\rom.

\rom{3)} The vacuum $\Omega=(1,0,0,\dots)\in\Fin(S)$ is a strong cyclic vector
of the family $(a^\sim(\xi))_{\xi\in S}$\rom, i\rom.e\rom{.,} the linear span of
the set
$$\{\Omega\}\cup\{\,a(\xi_1)\cdots a(\xi_n)\Omega\mid \xi_i\in S,\ i=1,\dots,n,\
n\in \N\,\}$$
is dense in $\Fin(S)$\rom.
\end{lem}

The item 1 of Lemma 2 we have already proved, while the item~2 follows from the
definition of the operators $a(\xi)$ and the estimates \eqref{2.1}--\eqref{2.3}.
The item~3 is implied by the fact that we have the standard creation operators,
which allows us just to repeat the proof of the corresponding fact in the
Gaussian and Poisson cases.

Thus, having Lemma~2, we apply the projection spectral theorem, which gives us
the existence of a spectral measure $\mu$ of the family $(a^\sim(\xi))_{\xi\in S}$ as
a probability measure on $S'$. For $\mu$-almost every $\om\in S'$, there exists
a generalized joint eigenvector
$$P(\om)=(P^{(n)}(\om))_{n=0}^\infty\in\Fin^*(S)$$
satisfying
\begin{equation}\label{2.9}
\la\!\la P(\om),a(\xi)f\ra\!\ra=\la\om,\xi\ra\la\!\la P(\om),f\ra\!\ra,
\qquad f\in \Fin(S).\end{equation}
Then, the operator $I$ defined by
\begin{gather}
If=(If)(\om)=\la\!\la P(\om),f\ra\!\ra=\sum_{n=0}^\infty\la
P^{(n)}(\om),f^{(n)}\ra\,n!,\label{2.10}\\
f=(f^{(n)})_{n=0}^\infty\in\Fin(S),\notag\end{gather}
can be extended by continuity to a unitary operator between $\Gama(\Ha)$ and\linebreak
$L^2(S',{\cal B}(S'),d\mu)$, under which any operator $a^\sim(\xi)$ goes over
into the operator of multiplication by $\la\om,\xi\ra$ (see \eqref{2.9} and
\eqref{2.10}). Denoting $\wom{n}=P^{(n)}(\om)\,n!$, we present the unitary $I$ in
the form \eqref{2.5}. The recurrence formula \eqref{2.6} can now be derived from
the equality
\begin{gather*}
\la\om,\xi\ra\la\wom{n},\xi^{\ot n}\ra=\la\wom{(n+1)},\xi^{\ot(n+1)}\ra
+\la\wom n,2n(\xi^2)\ho\xi^{\ot(n-1)}+\la\xi\ra\xi^{\ot n}\ra\\
\text{}+\la\wom{(n-1)},n\la\xi^2\ra\xi^{\ot(n-1)}+n(n-1)(\xi^3)\ho\xi^{\ot(n-2)}\ra.
\end{gather*}
As easily seen, $\wom{n}$ is well defined as an element of $S^{\prime\,\ho n}$
for {\it each\/} $\om\in S'$, not only for $\mu$-almost all $\om\in S'$.

Thus, we need only to prove that $\mu$ is, in fact, the measure of
Gamma white noise.

Let ${\cal B}_{\mathrm b}(\R^d)$ denote the set of all bounded Borel
sets in $\R^d$, let $\Delta\in{\cal B}_{\mathrm b}(\R^d)$, and let
$\chi_\Delta$ denote the indicator of $\Delta$. Denote by $a^\sim(\chi_\Delta)$
the operator in $\Gama(\Ha)$ whose image under the unitary $I$ is the operator of
multiplication by the function
$\la\om,\chi_\Delta\ra=\la\wom1,\chi_\Delta\ra+\sigma(\Delta)
=(I\chi_\Delta)(\om)+\sigma(\Delta)$.

Each of the vectors $\chi_\Delta^{\otimes n}$, $n\ge0$,
$\chi_\Delta^{\otimes 0}:=\Omega$,
belongs to $\Gama(\Ha)$, and let ${\cal K}_\Delta$ be the subspace of $\Gama(\Ha)$
spanned by these vectors. A limiting procedure shows that
\begin{equation}\label{2.12}a^\sim(\chi_\Delta)\chi_\Delta^{\ot
    n}=\chi_\Delta^{\ot(n+1)}+
(2n+\sigma(\Delta))\chi_\Delta^{\ot
  n}+n(n-1+\sigma(\Delta))\chi_\Delta^{\ot(n-1)}.\end{equation}
Therefore, ${\cal K}_\Delta$ is an invariant subspace for the operator
$a^\sim(\chi_\Delta)$.  Let $a_\Delta$ stand for the restriction of $a^\sim(\chi_\Delta)$ to the linear span of the vectors
$\chi_\Delta^{\otimes n}$. Then, $a_\Delta$ is a densely defined, Hermitian operator
in ${\cal K}_\Delta$.

Let $c_{n,\Delta}=\|\chi_\Delta^{\ot n}\|_{\Gama(\Ha)}$, then the vectors
$(e_\Delta^{(n)})_{n=0}^\infty$, with
$e_\Delta^{(n)}=c_{n,\Delta}^{-1}\chi_\Delta^{\ot n}$, form an orthonormal basis
(ONB) in ${\cal K}_\Delta$.
\eqref{2.12} yields
\begin{equation}
a_\Delta e^{(n)}_\Delta=\frac{c_{n+1,\Delta}}{c_{n,\Delta}}\,e_\Delta^{(n+1)}+
(2n+\sigma(\Delta))e_\Delta^{(n)}
+n(n-1+\sigma(\Delta))\,\frac{c_{n-1,\Delta}}{c_{n,\Delta}}\,e_\Delta^{(n-1)}.\label{2.13}
\end{equation}
Since $a_\Delta$ is Hermitian, the matrix of this operator in the ONB
$(e^{(n)}_\Delta)_{n=0}^\infty$ must be symmetric, which together with the
formula \eqref{2.13} implies that this matrix is a Jacobi one, see \cite{Be}, i.e.,
\begin{equation}\label{2.14}
a_\Delta e_\Delta^{(n)}=\alpha_{n+1,\Delta}e_\Delta^{(n+1)}+
\beta_{n,\Delta}e_\Delta^{(n)}+\alpha_{n,\Delta}
e_\Delta^{(n-1)}.
\end{equation}
Comparing \eqref{2.13} and \eqref{2.14} gives the equality
$$\frac{c_{n,\Delta}}{c_{n-1,\Delta}}=n(n-1+\sigma(\Delta))\,\frac{c_{n-1,\Delta}}{c_{n,\Delta}},$$
from where
\begin{equation}\label{2.15}
\frac{c_{n,\Delta}}{c_{n-1,\Delta}}=\sqrt{n(n-1+\sigma(\Delta))}.\end{equation}
Substituting \eqref{2.15} into \eqref{2.13}, we derive that the coefficients in
\eqref{2.14} equal
\begin{equation}\label{2.16}
\alpha_{n,\Delta}=\sqrt{n(n-1+\sigma(\Delta))},\qquad \beta_{n,\Delta}=2n+\sigma(\Delta).
\end{equation}

By using e.g.\ \cite{Be}, Ch.~7, Th.~1.3, we conclude that the Jacobi matrix
defined by the
coefficients \eqref{2.16} has a unique spectral measure. This means that there
exists a unique probability measure $\mu_\Delta$ on $\R$ such that the system of
polynomials $(P_{n,\Delta})_{n=0}^\infty$ defined by the recurrence formula
\begin{equation}\label{2.17}
P_{n+1,\Delta}(s)=\frac{1}{\alpha_{n+1,\Delta}}\big((s-\beta_{n,\Delta})P_{n,\Delta}(s)
-\alpha_{n,\Delta}P_{n-1,\Delta}(s)
\big)\end{equation}
forms an ONB in  $L^2(\R,d\mu_\Delta)$, and under the unitary $I_\Delta\colon
{\cal K}_\Delta\to L^2(\R,d\mu_\Delta)$ given by
$I_\Delta e_\Delta^{(n)}=P_{n,\Delta}$ the
operator $a_\Delta^\sim$---the closure in ${\cal K}_\Delta$ of the
essentially selfadjoint operator $a_\Delta$---goes
over into the operator of multiplication by
$s$. Moreover, the form of the coefficients \eqref{2.16} gives (see, e.g., \cite{Chihara}) that the measure
$\mu_\Delta$ is exactly the Gamma distribution having density \eqref{gamma}.

On the other hand, since $Ia^\sim(\chi_\Delta)I^{-1}$ is the operator of
multiplication by $\la\om,\chi_\Delta\ra$, \eqref{2.14} and \eqref{2.17} imply that
$$(Ie_\Delta^{(n)})(\om)=P_{n,\Delta}(\la\om,\chi_\Delta\ra).$$
Therefore, $\mu_\Delta$ coincides with the distribution of the r.v.\
$\la\om,\chi_\Delta\ra.$

Next, for arbitrary disjoint  sets $\Delta_1,\Delta_2\in{\cal B}
_{\mathrm b}(\R^d)$,
the spaces ${\cal K}_{\Delta_1}\ominus \widehat \Omega$ and ${\cal
  K}_{\Delta_2}\ominus\widehat\Omega$
are orthogonal in
$\Gama(\Ha)$. Here, $\widehat \Omega$  denotes the vacuum space spanned by
$\Omega$, and $\ominus$ denotes the orthogonal difference. Therefore,
the r.v.'s $\la\om,\chi_{\Delta_1}\ra$ and
$\la\om,\chi_{\Delta_2}\ra$ are independent.

Thus, $\mu$ is indeed the measure of Gamma white noise.\quad $\Box$\vspace{2mm}

\noindent {\it Remark\/} 5. Note that the polynomial $P_{n,\Delta}(s)$ is up to
a multiplier $(-1)^n$ the normalized Laguerre polynomial with parameter
$\sigma(\Delta)-1$.
Therefore,
$$\la\wom{n},f^{(n)}\ra=(-1)^n\la L_n^\sigma(\om),f^{(n)}\ra,$$
where $L_n^\sigma(\om)$ is the kernel of the generalized Laguerre polynomial
introduced in \cite{Silva}. This yields, in particular,
that the Gamma-Wick exponential defined by
\begin{equation}
\label{exponential}
\wexp{\om}{\varphi}=\sum_{n=0}^\infty\frac1{n!}\la\wom{n},\varphi^{\ot n}\ra,\qquad \varphi\in
S_\C,\end{equation}
has
the following
representation,
for $\varphi\in S$ such that $\varphi>-1$,
\begin{equation}\label{silva}
\wexp{\om}{\varphi}=\exp\Big[\la\om,\frac{\varphi}{\varphi+1}\ra-\la\log(1+\varphi)\ra\Big].\end{equation}
 \vspace{2mm}

\noindent {\it Remark\/} 6. Following the tradition of quantum probability,
e.g., \cite{Mey}, the family of the operators $(a^\sim(\chi_\Delta))_{\Delta\in\Lambda}$
can be called a quantum Gamma process in the extended Fock space $\Gama(\Ha)$.

\section{Chaos decomposition of the Gamma space}

In this section, we will discuss some properties of the unitary $I$ and the
space of Gamma white noise $\ltwoG$.

Denote by $\Pe(S')$ the set of continuous polynomials on $S'$, i.e., functions
on $S'$ of the form
\begin{equation*}
\vvv(\om)=\sum_{i=0}^n\la\om^{\ot i},f^{(i)}\ra,\qquad f^{(i)}\in S_\C^{\ho i},\
\om^{\ot 0}=1,\ i\in\N_0.\end{equation*}
The greatest number $n$ for which $f^{(n)}\ne0$ is called the power of a
polynomial. Denote also by $\Pe_n(S')$ the set of continuous polynomials of
power $\le n$.

By using the isomorphism
\begin{equation}\label{3.1}
\Fin(S)\ni f=(f^{(n)})_{n=0}^\infty\mapsto\sum_{n=0}^\infty\la\om^{\ot
  n},f^{(n)}\ra\in\Pe(S'),
\end{equation}
one induces from $\Fin(S)$ a topology on $\Pe(S')$, which makes it a nuclear
space.

\begin{prop} We have
$$I(\Fin(S))=\Pe(S'),$$
and the topology on $\Pe(S')$ induced from $\Fin(S)$ by $I$ coincides with that
induced by \eqref{3.1}\rom.\end{prop}

\noindent {\it Proof}. Using the recurrence relation
\eqref{2.6}, one obtains by induction the inclusion $I(\Fin(S))\subset\Pe(S')$,
and moreover, the equality
\begin{equation}\label{3.2}
\la\wom{n},f^{(n)}\ra=\la\om^{\ot n},f^{(n)}\ra+p_{n-1}(\om),\end{equation}
where $p_{n-1}(\om)\in\Pe_{n-1}(S')$. Using \eqref{3.2}, one obtains again by
induction the inverse inclusion ${\cal P}(S')\subset\Fin(S)$. The statement about
the topology on $\Pe(S')$ follows directly from the above
consideration. \quad$\Box$\vspace{2mm}

For every $n\in\N_0$, put $(L^2_{{\mathrm G},n})=I(\Gama^{(n)}(\Ha))$, so that
we obtain the following decomposition of the space $\ltwoG$:
\begin{equation}\label{3.3}
\ltwoG=\bigoplus_{n=0}^\infty(L^2_{{\mathrm G},n}).\end{equation}

\begin{prop}For each $n\in\N$\rom, we have
\begin{equation}\label{gamman}
(L^2_{{\mathrm G},n})=\Pe_{{\mathrm G},n}(S')\ominus\Pe_{{\mathrm G},n-1}(S'),\end{equation}
where $\Pe_{{\mathrm G},n}(S')$ denotes the closure of $\Pe_n(S')$ in the $\ltwoG$
norm\rom, and $\ominus$ stands for the orthogonal difference in $\ltwoG$.\end{prop}

\noindent {\it Proof}. As appears from the proof of Proposition~1,
$$I\Big(\bigoplus_{i=0}^n S_\C^{\ho i}\Big)=\Pe_n(S'),$$
and hence
$$I\Big(\bigoplus_{i=0}^n\Gama^{(n)}(\Ha)\Big)=\bigoplus_{i=0}^n(L^2_{{\mathrm
    G},n})=\Pe_{{\mathrm G},n}(S'),$$
from where \eqref{gamman} follows.\quad$\Box$

\begin{prop}
Let $P_{{\mathrm G},n}$ stand for the orthogonal projection of $\ltwoG$ onto the
subspace $(L_{{\mathrm G},n}^2)$\rom. Then
\begin{equation}\label{3.4}
P_{{\mathrm G},n}(\la\om^{\ot
  n},f^{(n)}\ra)=\la\wom{n},f^{(n)}\ra,\qquad f^{(n)}\in S_\C^{\ho n}.\end{equation}\end{prop}

\noindent{\it Proof}. \eqref{3.4} is derived from \eqref{3.2}
and Proposition~2
by
applying to the left and right hand sides of the
latter formula the operator $P_{{\mathrm G},n}$ and taking to notice that
$\la\wom{n},f^{(n)}\ra\in(L^2_{{\mathrm G},n})$.\quad$\Box$\vspace{2mm}

Propositions 2 and 3  make it possible to interpret \eqref{3.3} as a kind of  chaos
decomposition of $\ltwoG$. However, one should be  careful with this term.
Indeed, since ${\Bbb E}[X_\Delta]=\sigma(\Delta)$ and ${\Bbb E}[(X_\Delta-\sigma(\Delta))^2]=\sigma(\Delta)$,
 the compensated Gamma process
$$\tilde X_\Delta=\tilde X_\Delta(\om)=
X_\Delta(\om)-\sigma(\Delta)=\la\om,\chi_\Delta\ra-\sigma(\Delta)
=\la\wom1,\chi_\Delta
\ra
,\qquad \Delta\in{\cal B}_{\mathrm b}(\R^d),$$
is a normal martingale, and so one can introduce  multiple stochastic integrals
w.r.t.\ $\tilde X_\Delta$. But it follows from the general result of \cite{Der}
that this process does not possess the chaotic representation property, so that
\eqref{3.3} cannot be a chaos decomposition of $\ltwoG$ in the sense of multiple
stochastic integrals. To be more precise, we present the following propopsition.

\begin{prop}
The chaos in $\ltwoG$ generated by multiple stochastic integrals
w\rom.r\rom.t\rom.\ the compensated Gamma process $\tilde X_\Delta$
coincides with the image under the unitary $I$ of the usual Fock space ${\cal
  F}(\Ha)$ as a subspace of $\Gama(\Ha)$\rom, and moreover
$$ \int_{\R^{dn}}f^{(n)}(x_1,\dots,x_n)\,d\tilde X_{x_1}\cdots d\tilde
X_{x_n}=If^{(n)}$$
for an arbitrary $f^{(n)}\in\Ha_\C^{\ho n}$\rom, where $\Ha_\C^{\ho n}$ is
considered as a subspace of $\Gama^{(n)}(\Ha)$\rom.\end{prop}

\noindent{Proof}. We remind that the $n$-fold stochastic integral w.r.t.\ $\tilde
X_\Delta$ are defined as follows (e.g., \cite{Mey,Der}). First, one takes
arbitrary disjoint sets $\Delta_1,\dots,\Delta_n\in\Lambda$ and sets
\begin{equation}\label{3.5}
\int_{\R^{dn}}\big(
\chi_{\Delta_1}(x_1)\cdots\chi_{\Delta_n}(x_n)
\big)^\sim
d\tilde X_{x_1}\cdots d\tilde X_{x_n}
=\tilde X_{\Delta_1}\cdots \tilde X_{\Delta_n}.
\end{equation}
Then, one extends the equality \eqref{3.5} by linearity and continuity in
the $\Ha_\C^{\ho n}$ norm. But Theorem~1 yields that
\begin{gather*}
\tilde X_{\Delta_1}(\om)\cdots\tilde
X_{\Delta_n}(\om)=\la\wom1,\chi_{\Delta_1}\ra\cdots\la\wom{1},\chi_{\Delta_n}\ra\\
=I\big(
(a(\chi_{\Delta_1})-\sigma(\Delta_1)\id)\cdots(a(\chi_{\Delta_n})-\sigma(\Delta_n)\id)\Omega
\big)\\
=I(a^+(\chi_{\Delta_1})\cdots a^+(\chi_{\Delta_n})\Omega)
=I(\chi_{\Delta_1}\ho\cdots\ho\chi_{\Delta_n})
,\end{gather*}
which together with our consideration at the end of Section~1 gives the statement.\quad$\Box$

\section{Spaces of test and generalized functions\\
 and coordinate operators on
  them}
In this section, we will construct spaces of test and generalized functions of
Gamma white noise and consider some their properties.
 The simplest ones can be obtained by applying $I$ to \eqref{a}.
Thus, taking to notice Proposition~1, we get
$$\ltwoG\supset\Pe(S'),\qquad \Pe(S')^*\supset\Pe(S').$$
The dual $\Pe(S')^*$ of the space of continuous polynomials $\Pe(S')$ consists
of generalized functions, which we present in the form
$$\Phi=\Phi(\om)=\sum_{n=0}^\infty\la\wom n,F^{(n)}\ra,$$
where $F=(F^{(n)})_{n=0}^\infty\in\Fin^*(S)$ and the dualization of $\Phi$ with
$\vvv\in\Pe(S')$, $\vvv(\om)=\sum_{n=0}^\infty \la\wom n,f^{(n)}\ra$, is
given by
$$\la\!\la\Phi,\vvv\ra\!\ra =\sum_{n=0}^\infty\la \overline{F^{(n)}},f^{(n)}\ra n!.$$

However, the test space $\Pe(S')$ is too small. This is why we will consider
the following nuclear space and its dual: (cf.~\cite{KoLeS,Koetal,KSWY,GAS})
\begin{gather}
\Gama(\Ha)
\supset\projlim_{p,k\to\infty}{\cal F}_{1,k}(S_p)={\cal F}_1(S_1),\notag\\
{\cal F}_{-1}(S')=\indlim_{p,k\to\infty}{\cal F}_{-1,-k}(S_{-p})\supset
\projlim_{p,k\to\infty}{\cal F}_{1,k}(S_p)={\cal F}_1(S_1),
\label{4.1}\end{gather}
where the spaces ${\cal F}_{1,k}(S_p)$, $p\ge1$, $k\ge k_0$, are defined by \eqref{space}, ${\cal
  F}_{-1,-k}(S_{-p})$ are their respective duals (w.r.t.\ ${\cal F}(\Ha)$):
$${\cal F}_{-1,-k}(S_{-p})=\bigoplus_{n=0}^\infty S_{-p,\C}^{\ho n} 2^{-nk}.$$
The application of $I$
 to \eqref{4.1} gives the following spaces of test and generalized functions:
\begin{gather*}
\ltwoG\supset\projlim_{p,k\to\infty}(S_{\mathrm
  G})_{p,k}^1=(S_{\mathrm G})^1,\\
(S_{\mathrm G})^{-1}=\indlim_{p,k\to\infty}(S_{\mathrm G})^{-1}_{-p,-k}\supset
\projlim_{p,k\to\infty}(S_{\mathrm G})^1_{p,k}=(S_{\mathrm G})^1.
\end{gather*}

On the space $\SG^{-1}$, we introduce an $\S$-transform in a standard way
(cf.~\cite{KoLeS,Koetal,KSWY,GAS}): If $\Phi\in\SG^{-1}$, then by the definition of $\SG^{-1}$ there are
$p\ge1$ and $k\ge k_0$ such that $\Phi\in\SG_{-p,-k}^{-1}$, and we put
\begin{equation}\label{4.2}
 \S[\Phi](\theta)=\ll\overline{\Phi},\wexp{\cdot}{\theta}\rr,\qquad \theta\in
 S_\C,\ |\theta|_p<2^{-k/2},\end{equation}
where $\overline{\Phi}$ is the complex conjugate of $\Phi$ and  $\wexp{\cdot}{\theta}$
is defined by \eqref{exponential}.
The condition on the norm of $\xi$ in \eqref{4.2} implies that $\wexp{\cdot}\theta\in\SG_{p,k}^1
$, so that the dualization in \eqref{4.2} is well-defined. From the definition
  of the $\S$-transform, we have for each
$$\Phi=\Phi(\om)=\sum_{n=0}^\infty\la\wom{n},F^{(n)}\ra,$$
that
$$\S[\Phi](\theta)=\sum_{n=0}^\infty\la F^{(n)},\theta^{\ot
  n}\ra.$$

Notice that if $\theta\in S_\C'$, we can still  define the Wick exponential
$\wexp{\cdot}{\theta}$ by the same formula \eqref{exponential}, considering it
as a generalized function from $\SG^{-1}$. Hence, the $\S$-transform  of a test
function $\phi(\om)=\sum_{n=0}^\infty\la \wom{n},f^{(n)}\ra\in\SG^1$ can be
extended on $S_\C'$ by setting
\[\S[\phi](\theta)=\ll \overline{\phi},\wexp{\cdot}{\theta}\rr=\sum_{n=0}^\infty
\la f^{(n)},\theta^{\ot n}\ra,\qquad \theta\in S_\C'.\]

The following two theorems, which are due to \cite{KoLeS}, give the
description of the spaces $\SG^1$ and $\SG^{-1}$ in terms of the $\S$-transform.

\begin{th} For any $\Phi\in\SG^{-1}$\rom, $\S[\Phi]$ is a function holomorphic
  at zero\rom, i\rom.e\rom{.,} $\S[\Phi]\in\Hol$\rom. And conversely\rom, for any
  $F\in\Hol$\rom, there exists a unique $\Phi\in\SG^{-1}$ such that
  $\S[\Phi]=F$ \rom(i\rom.e\rom{.,} there is a neighborhood of zero\rom, $U$\rom,
  such that $\S[\Phi](\theta)=F(\theta)$ for all $\theta\in U$\rom{).}
\end{th}

Denote by $\Emin{S_\C'}$ the set of entire functions on $S_\C'$ of the first
order of growth and of minimal type, i.e., a function $u$ entire on $ S_\C'$
belongs to $\Emin{S_\C'}$ iff
\[ \forall p>0\ \forall \epsilon>0\ \exists C>0:\quad |u(\theta)|\le
Ce^{\epsilon|\theta|_{-p}}, \qquad \theta\in S_{-p,\C}.\]

The set $\Emin{S_\C'}$ can be topologized by the following family of norms
\begin{equation}\label{complexnorm}
\chush u\chush_{1,p,k}=\sup_{\theta\in S_{-p,\C}}\Big\{|u(\theta)|\exp\big[-\frac
1k\,|\theta|_{-p}\big]\Big\},\qquad p,k\in\N,\end{equation}
so that $\Emin{S_\C'}$ becomes a countably normed space.

\begin{th}
For any $\vvv\in\SG^{1}$\rom, $\S[\vvv]\in\Emin{S_\C'}$\rom. And
conversely\rom, for any $f\in\Emin{S_\C'}$\rom, there exits a unique $\vvv\in
\SG^1$ such that $\S[\vvv]=f$\rom. Moreover\rom, the $\S$-transform is a
homeomorphism between the topological spaces $\SG^1$ and $\Emin{\S_\C'}.$
\end{th}

Denote now by $\Emin{S'}$ the set of the restrictions to $S'$
 of  functions from $\Emin{S_\C'}$.
The topology on $\Emin{S'}$ is induced by that of $\Emin{S'_\C}$.

The next theorem gives the inner description of the test space $\SG^1$.

\begin{th}
We have
\begin{equation}\label{4.3} \SG^1=\Emin{S'},\end{equation}
where \eqref{4.3} is understood as an equality of topological
spaces\rom.\end{th}

\noindent{\it Proof}.
This theorem is, in fact, a direct corollary of results of  \cite{GAS} (see also
\cite{LyUs}),
since the Gamma-Wick monomials
$\la\wom{n},f^{(n)}\ra$
 form a generalized Appell system in the sense of \cite{GAS} with the
 transformation function $\alpha(\varphi)=\dfrac{\varphi}{\varphi+1}$
(cf.~\cite{Silva}).
However, for the further exposition, we need to state some details of the theory
of (generalized) Appell polynomials.

Let $\mu$ be a probability measure on $(S',{\cal B}(S'))$ satisfying the following
two conditions:\vspace{1.7mm}

\noindent{\bf Condition 1.} The Laplace transform $\ell_\mu(\varphi)$
of the measure $\mu$ can be extended to an analytic function in a neighborhood
of zero in $S_\C'$. , \vspace{1.7mm}

\noindent{\bf Condition 2.} If $\phi\in{\cal P}(S')$ and $\phi=0$ $\mu$-almost
everywhere, then $\phi\equiv 0$.\vspace{1.7mm}

Put
$$e_\mu(\ttt;\om)=\frac{e^{\la\om,\ttt\ra}}{\ell_\mu(\ttt)}.$$
Due to Condition 1 and the fact that $\ell_\mu(0)=1$, for every fixed $\om\in
S'$
$e_\mu(\varphi;\om)$ is
an analytic function of $\ttt$ in a neighborhood of zero.

Define now $P^{(0)}_\mu(\om)=1$ and $P_\mu^{(n)}(\om)\in S^{\prime\,\ho n}$,
$n\in \N$, by
\begin{equation}\label{zyp}
\la P_\mu^{(n)}(\om),\ttt^{\ot
  n}\ra =\frac{d^n}{dt^n}\Big|_{t=0}e_\mu(t\varphi;\om),\end{equation}
so that
$$e_\mu(\varphi;\om)=\sum_{n=0}^\infty\frac{1}{n!}\,\la
P^{(n)}_\mu(\om),\ttt^{\ot n}\ra.$$
By using the kernel theorem, one shows that this definition is correct, i.e.,
that \eqref{zyp} really determines an element of $S^{\prime\,\ho n}$. A function
of $\om$ of the form $\la P_\mu^{(n)}(\om),f^{(n)}\ra$, where $f^{(n)}\in S_\C^{\ho
  n}$, is called an Appell polynomial. These polynomials have the following very important property:
let $\upsilon\in S'$ and let $D_\upsilon$ denote the G\^ateaux derivative in
direction $\upsilon$:
$$D_\upsilon \phi(\om)
=\frac{d}{dt}\Big|_{t=0}\,\phi(\om+t\upsilon)
=\lim_{t\to0}\frac{\phi(\om+t\upsilon)-\phi(\om)}{t},$$
then
$$D_\upsilon \la P^{(n)}_\mu(\om),\ttt^{\ot n}\ra=n\la\upsilon,\ttt\ra\la
P^{(n-1)}_\mu(\om), \ttt^{\ot(n-1)}\ra,$$
i.e., $D_\upsilon$ acts as an annihilation operator w.r.t.\ the Appell
polynomials. Particularly, for the gradient
\begin{equation}\label{gradient}
\nabla_x=D_{\delta_x},\qquad x\in\R^d,\end{equation}
we have
$$\nabla_x\la P_\mu^{(n)}(\om),\ttt^{\ot n}\ra=n\ttt(x)\la
P^{(n-1)}_\mu(\om),\ttt^{\ot(n-1)}\ra.$$
and so
$$\nabla_x e_\mu(\ttt;\om)=\ttt(x)e_\mu(\ttt;\om).$$

Now, for each $p\ge1$ and $k\in\N_0$, one introduces a norm $\norm\cdot\norm_{1,p,k}$ on ${\cal P}(S')$
as follows: for
$$\phi(\om)=\sum_{n=0}^k\la P_\mu^{(n)}(\om),f^{(n)}\ra,$$
we have
$$\norm\phi\norm_{1,p,k}^2=\sum_{n=0}^\infty |f^{(n)}|_p^2(n!)^2 2^{nk},$$
and let $[S_\mu]_{p,k}^1$
denote the closure of ${\cal P}(S')$ in this norm. Put
$$[S_\mu]^1=\projlim_{p,k\to\infty}[S_\mu]_{p,k}^1.$$
 Thus, one has \cite{KSWY} that there exist $p_1\ge1$ and $k_1\in\N_0$ such that $[S_\mu]^1_{p_1,k_1}$
 is topologically embedded into $L^2(S',d\mu)=(L^2_\mu)$, and moreover, one has
 the topological equality
$$[S_\mu]^1={\cal E}^1_{\mathrm min}(S').$$

Next, we proceed to consider a generalized Appell system. Let $\alpha$ be an
arbitrary function that maps a neighborhood of zero in $S_\C$ into itself, which
is supposed to be holomorphic, invertible and to satisfy $\alpha(0)=0$.

Let
$$e_{\mu,\alpha}(\ttt;\om)=e_\mu(\alpha(\ttt);\om)=\frac{e^{\la
    \om,\alpha(\ttt)}\ra}{\ell_\mu(\alpha(\ttt))}.$$
Then, analogously to the above, one defines generalized Appell polynomials\linebreak $\la
    P^{(n)}_{\mu,\alpha}(\om),f^{(n)}\ra$  in such a way that
$$e_{\mu,\alpha}(\ttt;\om)=\sum_{n=0}^\infty \frac1{n!}\,\la
P^{(n)}_{\mu,\alpha}(\om),\ttt^{\ot n}\ra.$$

The $\alpha$-gradient $\nabla_x^\alpha$ is defined now by
$$\nabla_x^\alpha \la P^{(n)}_{\mu,\alpha}(\om),\ttt^{\ot n}\ra=n\ttt(x)\la
P^{(n-1)}
_{\mu,\alpha}
(\om),\ttt^{\ot(n-1)}\ra,$$
so that
$$\nabla_x^\alpha e_{\mu,\alpha}(\ttt;\om)=\ttt(x)e_{\mu,\alpha}(\ttt;\om).$$

Using the polynomials $\la P_{\mu,\alpha}^{(n)}(\om),f^{(n)}\ra$, one defines
just as above the spaces $[S_\mu]_{p,k,\alpha}^1$, and again \cite{GAS} there
exist $p_2\ge 1$ and $k_2\in\N_0$ such that $[S_\mu]^1_{p_2,k_2,\alpha}$ is topologically embedded
into $(L^2_\mu)$ and
$$[S_\mu]_\alpha
^1=\projlim_{p,k\to\infty}[S_\mu]^1_{p,k,\alpha}={\cal E}^1_{\mathrm
  min}(S').$$

Let us come back to the case of the Gamma measure $\mu=\muG$. Evidently the
Laplace transform $\ell_{\mathrm G}(\ttt)$ of $\muG$ satisfies Condition~1 (see
\eqref{2.4a}).  Proposition~1 implies that Condition~2 is also satisfied for
$\muG$. Thus, we can put
\begin{equation}\label{sp}e_{\mathrm G}(\ttt;\om)=\frac{e^{\la\om,\ttt}\ra}{\ell_{\mathrm
    G}(\ttt)}=\exp\big[ \la\om,\ttt\ra+\la\log(1-\ttt)\ra\big].\end{equation}
By setting
\begin{equation}\label{spp} \alpha(\ttt)=\frac{\ttt}{\ttt+1},\end{equation}
which satisfies the conditions on $\alpha$, we get
\begin{equation}\label{neznaju}
e_{{\mathrm
    G},\alpha}(\ttt;\om)=\exp\Big[\big\la\om,\frac{\ttt}{\ttt+1}\big\ra-\la\log(1+\ttt)\ra\Big].\end{equation}
Comparing the last formula with \eqref{silva}, we conclude that
$$\la\wom{n},f^{(n)}\ra=\la P^{(n)}_{{\mathrm G},\alpha}(\om),f^{(n)}\ra,$$
which yields the theorem.
\quad $\blacksquare$\vspace{2mm}

On the space $\SG^{-1}$, one can introduce   a Wick product $\ds$ as follows
(cf.~\cite{KoLeS,KSWY,GAS}):
for arbitrary $\Phi,\Psi\in\SG^{-1}$, $\Phi\,\ds\,\Psi$ is an element of $\SG^{-1}$
such that
\begin{equation}\label{4.4}
  \S[\Phi\,\ds\,\Psi](\theta)=\S[\Phi](\theta)\,\S[\Psi](\theta)\end{equation}
for $\theta\in S_\C$ from some neighborhood of zero. Since $\Hol$
is an algebra under pointwise multiplication of functions, this definition is
correct. \eqref{4.4} implies that
\[\la\wom{n},F^{(n)}\ra\,\ds\,\la\wom{m},G^{(m)}\ra=\la\wom{(n+m)},F^{(n)}\ho
G^{(m)}\ra.\]
Notice that the Wick product of two test functions from $(S_{\mathrm G})^1$
belongs again to $(S_{\mathrm G})^1$.
\vspace{2mm}

\noindent{\it Remark\/} 7. Because of \eqref{sp}--\eqref{neznaju}, the above introduced Wick
product on $\SG^{-1}$ coincides with the Wick product on this space defined in
the framework of biorthogonal analysis \cite{KSWY}, i.e., by using the
function $e_{\mathrm G}(\cdot\,;\theta)$ instead of $\wexp{\cdot}{\theta}$.\vspace{2mm}

Let us consider a pair of simplest examples of generalized functions from $\SG^{-1}$.
The first one is the Gamma white noise:
\[
\om(x)={:}\,\om\,{:}_{\mathrm
  G}(x)+1=\la\wom{1},\delta_x\ra+1\in\SG_{-1,-k_0}^{-1},\qquad x\in\supp\sigma,
\]
the function $\om(x)-1={:}\,\om\,{:}_{\mathrm G}(x)$ can be thought of as the
compensated Gamma white noise.
Next, by taking the Wick product of ${:}\,\om\,{:}_{\mathrm G}(x_i)$, $i=1,\dots,n$,
we obtain a Gamma white noise monomial as
\begin{gather*}
{:}\,\om\,{:}_{\mathrm G}(x_1)\,\ds\,
{:}\,\om\,{:}_{\mathrm G}(x_2)\,
\ds\cdots\ds\,
{:}\,\om\,{:}_{\mathrm G}(x_n)=\la\wom{n},\delta_{x_1}\ho\cdots\ho\delta_{x_n}\ra\\
={:}\,\om^{\ot n}\,{:}_{\mathrm G}(x_1,\dots,x_n)\in\SG_{-1,-k_0}^{-1}.
\end{gather*}

The delta function of Gamma white noise, $\tilde\delta_{\upsilon}$, where
$\upsilon\in S'$, is defined by
\[
\ll\tilde\delta_\upsilon,\vvv\rr=\vvv(\upsilon),\qquad
\vvv\in\SG^{1}.\]
Evidently, $\tilde\delta_\upsilon$ belongs to $\SG^{-1}$ and has
the representation:
\[\tilde \delta_\upsilon=\tilde\delta_\upsilon(\om)=\sum_{n=0}^\infty \la\wom{n},\frac1{n!}\,\wup{n}
\ra.\]

Now, we wish to introduce operators of coordinate multiplication,
$\om(x){\cdot}$, acting from $\SG^1$ into $\SG^{-1}$.
To this end, we define first linear operators $\di_x\colon\SG^1\to\SG^1$ and
$\dig_x\colon \SG^{-1}\to\SG^{-1}$ for each $x\in\supp\sigma$ by
\begin{align*}
\di_x\la\wom{n},f^{(n)}\ra&=n\la\wom{(n-1)},f^{(n)}(x,\cdot)\ra,\\
\dig_x\la\wom{n},F^{(n)}\ra&=\la\wom{(n+1)},\delta_x\ho F^{(n)}\ra.
\end{align*}
It easy to show that these operators are continuous.

Let us preserve the same notations $a^+(\xi)$, $a^0(\xi)$, $a_1^-(\xi)$, and
$a_2^-(\xi)$ for the images of the corresponding operators under the unitary
$I$. Then, the operators $a^+(\xi)$ and $a_1^-(\xi)$ have the integral
representation (cf.\ \cite{Hidaetal}):
\begin{align}
a_1^-(\xi)&=\int_{\supp\sigma}d\sigma(x)\,\xi(x)\di_x,\notag\\
a^+(\xi)&=\int_{\supp\sigma}\,d\sigma(x)\,\xi(x)\dig_x.\label{4.5}\end{align}

Integrals of such  type are understood usually in the sense that one applies
pointwisely the integrand operator to a test function, then dualizes the result
with another test function, and finally integrates the obtained function of $x$
w.r.t.\ the measure $\sigma$.

As well known, the neutral operator $a^0(\xi)$ has the representation
\[a^0(\xi)=\int_{\supp\sigma}d\sigma(x)\, \xi(x)\,\dig_x\di_x.\]

A new point appearing in Gamma analysis is, of course, the second annihilation
operator, $a_2^-(\xi)$, which has now the representation
\begin{equation}\label{4.7}
a_2^-(\xi)=\int_{\supp\sigma}d\sigma(x)\,\xi(x)\,\dig_x\di_x\di_x  .\end{equation}
The formula \eqref{4.7} can be verified in a standard way.

Thus, we define the operator $\om(x){\cdot}$ by
\begin{equation}
\om(x){\cdot}=\dig_x+2\dig_x\di_x+1+\di_x+\dig_x\di_x\di_x,\label{multiplication}\end{equation}
which acts continuously from $(S_{\mathrm G})^1$ into $(S_{\mathrm G})^{-1}$.
Evidently,
\[
\la\om,\xi\ra{\cdot}=
\int_{\supp\sigma}d\sigma(x)\,\xi(x)\om(x){\cdot}
\, .\]

By analogy with \cite{Huang} and \cite{LyReShch}, the family of operators
$(\om(x){\cdot})_{x\in\supp\sigma}$ can be also  called a Gamma field, or a
quantum Gamma white noise process.

As will be shown now, the action of all the above operators can be easily
represented in terms of the $\S$-transform. Indeed,
let $\phi\in\SG^1$ and let $U(\theta)=\S[\phi](\theta)$.
The formula for $\dig_x$, resp.\ $a^+(\xi)$ is
well known (e.g., \cite{HKPS}):
\begin{alignat}{2}
\S[\dig_x\vvv](\theta)&=\theta(x)U(\theta),&&\qquad \theta\in S_\C,\notag\\
\S[a^+(\xi)\phi](\theta)&=\la\xi,\theta\ra U(\theta),&&\qquad \theta\in S_\C'
\notag\end{alignat}
(of course, the formula for $\dig_x$ holds for $\Phi\in \SG^{-1}$ with $\theta$
from a neiborhood of zero in $S_\C$).

Next, analogously to \cite{ItKu}, we have that
\begin{alignat}{2}
\S[(\di_x+2\dig_x\di_x)\vvv](\theta)&=D_{\delta_x(1+2\theta)}U(\theta),&&\qquad
\theta\in S_\C,\label{eins}\\
\S[(a_1^-(\xi)+2a^0(\xi))\phi](\theta)&=D_{\xi(1+2\theta)}U(\theta),&&\qquad
\theta\in S_\C',\label{zwei}\end{alignat}
where $D_\upsilon$, $\upsilon\in S_\C'$, denotes also the G\^ateaux derivatives
in direction $\upsilon$ of a function defined on $S_\C'$. Indeed, let $\phi$ be
of the form $\phi(\om)=\la\wom{n},\varphi^{\ot n}\ra$, then
$U[\phi](\theta)=\la\varphi^{\ot n},\theta^{\ot n}\ra$, and
\begin{gather*}
\frac d{dt}\Big|_{t=0}\la\varphi^{\ot n},(\theta+t\delta_x(1+2\theta))^{\ot n}\ra
=n\varphi(x)(1+2\theta(x))\la \varphi^{\ot(n-1)},\theta^{\ot(n-1)}\ra\\
=\S[n\varphi(x)\la\wom{(n-1)},\varphi^{\ot(n-1)}\ra+2n\varphi(x)\la\wom
n,\delta_x\ho\varphi^{\ot(n-1)}\ra](\theta)\\
=\S[(\di_x+2\dig_x\di_x)\la\wom n,\varphi^{\ot n}\ra](\theta),\qquad \theta\in S_\C.
\end{gather*}

Hence, the formula \eqref{eins} will be proved if we show that, if
$\phi_m\to\phi$ in $\SG^1$ as $m\to\infty$ and $U_m(\theta)=\S[\phi_m](\theta)$,
then for any fixed $\theta\in S_\C$,
\begin{equation}\label{drei}
\lim_{m\to\infty}D_{\delta_x(1+2\theta)}U_m(\theta)=D_{\delta_x(1+2\theta)}U(\theta).\end{equation}
To this end, we will use the following lemma:
\begin{lem}
For each $\upsilon\in S_\C'$\rom, $D_\upsilon$ defines a linear continuous
operator on $\Emin{S_\C'}.$\end{lem}

\noindent {\it Proof}. Every function $u\in\Emin{S_\C'}$ can be
represented in the form (see e.g.\ \cite{KoLeS})
\[ u(\theta)=\sum_{n=0}^\infty \la f^{(n)},\theta^{\ot n}\ra,\qquad f^{(n)}\in
S_\C^{\ho n}.\]
For each $p\ge1$ and $k\in\N_0$, define the norm
\begin{equation}\label{drym}
\haha u\haha_{1,p,k}^2=\sum_{n=0}^\infty |f^{(n)}|_p^2(n!)^2 2^{nk},\end{equation}
and let ${\cal E}^1_{p,k}(S'_\C)$ be the closure of $\Emin{S_\C'}$ in this norm.
Theorem~3 (see also \cite{KoLeS}) implies that
\[\Emin{S_\C'}=\projlim_{p,k\to\infty}{\cal E}^1_{p,k}(S_\C'),\]
i.e., that the sequence of norms \eqref{drym} is equivalent to
\eqref{complexnorm}
Let $\upsilon\in S_{-p',\C}$ with $p'\ge1$. Since
\[ D_\upsilon \la\varphi^{\ot n},\theta^{\ot n}
\ra=n\la\varphi,\upsilon\ra\la\varphi^{\ot(n-1)},\theta^{\ot(n-1)}\ra,\]
the norm of $D_\upsilon$ on each $({\cal E})_{p,k}^1$, $p\ge p'$, $k\in\N_0$,
is not greater than $|\upsilon|_{-p}2^{-k/2}$, and therefore $D_\upsilon$ acts
continuously on $\Emin{S_\C'}$.\qquad $\blacksquare$\vspace{2mm}

Now, from  Theorem~3 we conclude that the convergence $\phi_m\to\phi$ in $\SG^1$
implies the convergence $U_m\to U$ in $\Emin{S_\C'}$, and hence by Lemma~3
$D_\upsilon U_m\to D_\upsilon U$ in $\Emin{S_\C'}$ for any fixed $\upsilon\in
S_\C'$. Therefore, by \eqref{complexnorm} $D_\upsilon U_m(\theta)\to
D_\upsilon U(\theta)$ for any fixed $\theta\in S_\C'$, which evidently implies
\eqref{drei} because $\delta_x(1+2\theta)\in S_\C'$ for any $\theta \in S_\C$.
The formula \eqref{zwei} can be proved absolutely analogously to
\eqref{eins} if we take into account that $\xi(1+2\theta)\in S'_\C$ for an
arbitrary $\theta\in S_\C'$.

Again analogously to \eqref{eins}, we have
\[ S[\dig\di_x\di_x\phi](\theta)=\theta(x)\nabla_x^2 U(\theta),\qquad \theta\in
S_\C,\]
where $\nabla_x$ is defined by \eqref{gradient}. Then
\begin{equation}\label{vier}
\S[a_2^-(\xi)\phi](\theta)=\la\xi(x)\nabla_x^2 U(\theta),\theta(x)\ra,\qquad \theta\in
S_\C',\end{equation}
where $x$ denotes the variable in which the dualization is carried out.

This formula is evidently true for $\phi(\om)=\la\wom n,\varphi^{\ot n}\ra$: then
$U(\theta)=\la\varphi^{\ot n},\theta^{\ot n}\ra$ and
\begin{align*}
 \la \xi(x)\nabla_x^2\la\varphi^{\ot n},\theta^{\ot n}\ra,\theta(x)\ra
&       =\la
\xi(x)n(n-1)\varphi(x)^2\la\varphi^{\ot(n-2)},\theta^{\ot(n-2)}\ra,\theta(x)\ra\\
&      =n(n-1)\la\varphi^{\ot(n-2)},\theta^{\ot(n-2)}\ra\la\xi\varphi^2,\theta\ra\\
&      =n(n-1)\la \varphi^{\ot(n-2)}\ho(\xi\varphi^2),\theta^{\ot(n-1)}\ra.
 \end{align*}
Hence, it suffices to show that, if $U_m\to U$ in $\Emin{S_\C'}$, then for any
fixed $\theta\in S_\C'$ $\nabla_x^2U_m(\theta)\to\nabla_x^2U(\theta)$ in each
$S_{p,\C}$, $p\ge1$, as a function of $x$.

Representing
\begin{equation}\label{fuenf}
U_m(\theta)=\sum_{n=0}^\infty \la f_m^{(n)},\theta^{\ot n}\ra\to
U(\theta)=\sum_{n=0}^\infty\la f^{(n)},\theta^{\ot n}\ra\quad\text{in
  }\Emin{S_\C'},\end{equation}
we get
\begin{align*}
\nabla_x^2 U_m(\theta)&=\sum_{n=0}^\infty \nabla_x^2\la f_m^{(n)},\theta^{\ot
  n}\ra\\
&=\sum_{n=2}^\infty n(n-1)\la f^{(n)}_m(x,x,\cdot),\theta^{\ot(n-2)}\ra.
\end{align*}
\eqref{fuenf} yields that $f_m^{(n)}\to f^{(n)}$ in $S_\C^{\ho n}$, and
therefore
\[ \la f_m^{(n)}(x,x,\cdot),\theta^{\ot(n-2)}\ra\to\la
f^{(n)}(x,x,\cdot),\theta^{\ot(n-2)}\ra\quad\text{in }S_\C.\]
Moreover, upon \eqref{1.3} we have, for an arbitrary $\upsilon\in S_{-p,\C}$,
\begin{align*}
&\big| \la \upsilon(x),\la f_m^{(n)}(x,x,\cdot),\theta^{\ot(n-2)}\ra\ra
\big|\\
&\qquad =\big|\la
f_m^{(n)}(x_1,x_1,x_2,\dots,x_{n-1}),\upsilon(x_1)\theta^{\ot(n-2)}(x_2,\dots,x_{n-1})\ra\big|\\
&\qquad\le |\upsilon|_{-p}|\theta|_{-p}^{n-2}C_p|f_m^{(n)}|_p,
\end{align*}
whence taking to notice that $S_{p,\C}$ can be thought of as the dual of
$S_{-p,\C}$, we get
\[
\big|\la f_m^{(n)}(x,x,\cdot),\theta^{\ot (n-2)}\ra\big|_p\le
C_p|\theta|_{-p}^{n-2}|f_m^{(n)}|_p.\]
Estimating
\begin{align*}
&\sum_{n=2}^\infty n(n-1)\big|\la
f_m^{(n)}(x,x,\cdot),\theta^{\ot(n-2)}\ra\big|_p\\
&\qquad \le C_p\sum_{n=2}^\infty n(n-1)|\theta|_{-p}^{n-2}|f_m^{(n)}|_p\\
&\qquad \le C_p\sum_{n=2}^\infty \big((n-2)!\big)^{-1}|\theta|_{-p}^{n-2}|f_m^{(n)}|_pn!\\
&\qquad \le  C_p\Big(\sum_{n=0}^\infty (n!)^{-2}|\theta|_{-p}^{2n}\Big)^{1/2}\haha
U_m\haha_{1,p,0},
\end{align*}
we obtain the desired statement.

Thus, we have proved the following theorem.

\begin{th}
Let $\phi\in\SG^1$\rom, then the action of the operators
$\om(x){\cdot}$\rom, $x\in\supp\sigma$\rom,  and $\la\om,\xi\ra\cdot$\rom, $\xi\in
S$\rom, can be represented in terms of the
$\S$-transform as follows\rom:
\begin{alignat*}{2}
\S[\om(x){\cdot}\,\vvv](\theta)&=(\theta(x)+1)U(\theta)+D_{\delta_x(1+2\theta)}U(\theta)
+\theta(x)\nabla_x^2 U(\theta),&\qquad \theta\in S_\C,
\\
\S[\la\om,\xi\ra\cdot\phi](\theta)&=\la\xi,\theta+1\ra
U(\theta)+D_{\xi(1+2\theta)}U(\theta)+
\la\xi(x)\nabla_x^2U(\theta),\theta(x)\ra,&\qquad
\theta\in S_\C',\end{alignat*}
where $U(\theta)=\S[\phi](\theta).$
\end{th}

\section{Standard annihilation operator\\ on Gamma space}

In this section, we will study the standard annihilation operators $\di_x$ and
$a_1^-(\xi)$.

First, we note that
\[ \sup_{x\in\supp\sigma}|\delta_x|_{-p}=\|\delta\|_{p,\infty}<\infty.\]
Therefore, for each $x\in\supp\sigma$, $\di_x$ can be extended to a continuous
operator on $\SG^1_{p,k}$, $p\ge1$, $k\ge k_0$, with norm $\le \|\delta\|_{p,\infty} 2^{-k/2}$.

It follows from the proof of Theorem~4 that $\di_x$ is nothing but the $\alpha$-gradient
 with $\alpha(\varphi)=\frac{\varphi}{\varphi+1}$, and on the total set in
 $\SG_{p,k}^1$ consisting of the Wick exponentials $\wexp{\cdot}{\varphi}$
 with $|\varphi|_p<2^{-k/2}$ we have
\[
 \di_x\wexp{\om}{\varphi}=\nabla_x^\alpha\wexp{\om}{\varphi}=\varphi(x)\wexp{\om}{\varphi}.\]
Note that
\[ \alpha^{-1}(\varphi)=\frac{\varphi}{1-\varphi}=\sum_{n=1}^\infty \varphi^n,\]
and hence in virtue of \eqref{1.3}
\[ |\alpha^{-1}(\varphi)|_p\le\frac{|\varphi|_p}{1-C_p|\varphi|_p}\quad\text{if
  }|\varphi|_p<C_p.\]
Choosing $\varphi\in S$ such that
\[ |\varphi|_p<\frac{2^{-k/2}}{1+C_p 2^{-k/2}},\]
we get that $|\alpha^{-1}(\varphi)|_p<2^{-k/2}$,
and then upon \eqref{sp}--\eqref{neznaju}
\begin{align*}
\di_x e_{\mathrm G}(\varphi;\om)&=\di_x\wexp{\om}{\alpha^{-1}(\varphi)}\\
&=\frac{\varphi(x)}{1-\varphi(x)}\, e_{\mathrm G}(\varphi;\om)\\
&=\sum_{n=1}^\infty \varphi(x)^n \, e_{\mathrm G}(\varphi;\om).\end{align*}
Since
\[ \alpha(\ttt)=\frac{\ttt}{\ttt+1}=\sum_{n=1}^\infty
 (-1)^{n+1}\ttt^n,\]
the functions $e_{\mathrm G}(\ttt;\om)$ also constitute a total set in
$\SG^1_{p,k}$, and therefore, at least formally, we can write down
\begin{equation}\label{5.1}
\di_x =\alpha^{-1}(\nabla_x)
=\frac{\nabla_x}{1-\nabla_x}
=\sum_{n=1}^\infty \nabla_x^n.\end{equation}

Let ${\cal E}^{1}_{p,k}(S')$ denote the Hilbert space constructed in the same
way as ${\cal E}^1_{p,k}(S_\C')$
(see the proof of Lemma~3)
but only  starting from $\Emin{S'}$. Choosing an arbitrary $\tilde p\ge1$ and
$\tilde k\in \N_0$ such that $2^{\tilde k/2}>\|\delta\|_{\tilde p,\infty}$, we
get that each $\nabla_x$ is a continuous operator on ${\cal E}^1_{\tilde
  p,\tilde k}(S')$ with norm less than one. Hence, the series $\sum_{n=1}^\infty
\nabla_x^n$ converges in operator norm. Let now $p\ge1$ and $k\ge k_0$ be
chosen so that the space $\SG_{p,k}^1$ is topologically embedded into ${\cal
  E}^1_{\tilde p, \tilde k}(S')$. Then, the above functions $e_{\mathrm
  G}(\ttt;\om)$ constitute also a total set in ${\cal E}^1_{\tilde p,\tilde
  k}(S')$, which implies the equality \eqref{5.1} on each space ${\cal
  E}^1_{\tilde p,\tilde k}(S')$, and therefore on $\Emin{S'}$.

In the same way, one can derive the inverse representation of \eqref{5.1}
\begin{equation}\label{5.2}
\nabla_x =\alpha(\di_x)=\frac{\di_x}{\di_x+1}=\sum_{n=1}^\infty
(-1)^{n+1}\di_x^n.\end{equation}
As a corollary of \eqref{5.1} or \eqref{5.2}, we have the commutation of
arbitrary $\nabla_{x_1}$ and $\di_{x_2}$ on $\SG^1$.

Now, we will show that analogously to the one-dimensional case (see e.g.\
\cite{OGF}) the operators $\di_x$ and $a_1^-(\xi)$ have a representation as an integral
w.r.t.\ a difference operator.

\begin{th}
For an arbitrary $\vvv\in \SG^1$, we have
\begin{align}
(\di_x\vvv)(\om)&=\int_0^\infty\big(\vvv(\omega+s\delta_x)-\vvv(\omega)\big)e^{-s}\,ds\notag\\
&=\int_0^\infty \vvv(\omega+s\delta_x)e^{-s}\,ds-\vvv(\om),\qquad
x\in\supp\sigma,\label{5.1a}\\
(a_1^-(\xi)\vvv)(\om)&=\int_{\supp\sigma}\int_0^\infty\xi(x)\big(\vvv(\omega+s\delta_x)-\vvv(\om)\big)
e^{-s}\,ds\,d\sigma(x),\qquad \xi\in S.\label{5.2a}
\end{align}
\end{th}

\noindent {\it Proof}.
Let $\theta\in S_{-p,\C}$, then for any $u\in\Emin{S_\C'}$,
\begin{align*}
\Big|\int_0^\infty
u(\theta+s\delta_x)e^{-s}\,ds\Big|&\le\chush u\chush_{1,p,k}\int_0^\infty\exp
\big[\tfrac1k|\theta+s\delta_x|_{-p} \big] e^{-s}\,ds
\\
&\le\chush u\chush_{1,p,k}\exp\big[\tfrac1k|\theta|_{-p}\big]\int_0^\infty
\exp\big[ -s\big(1-\tfrac1k
\|\delta|_{p,\infty}
\big)\big]\,ds\\
&=\chush u\chush_{1,p,k}
\exp\big[\tfrac1k|\theta|_{-p}\big]\,\frac1{1-\frac
  1k\|\delta\|_{p,\infty}},\qquad k>\|\delta\|_{p,\infty},
\end{align*}
where the norm $\chush\cdot\chush_{1,p,k}$ is defined by \eqref{complexnorm}.
Hence
\[\Chush\int_0^\infty u(\cdot+s\delta_x)e^{-s}\,ds\,\Chush_{1,p,k}\le
\chush u\chush_{1,p,k}\,
\frac{1}{1-\frac 1k\|\delta\|_{p,\infty}}.\]
Therefore, by Theorem~4 the operator $A_x$ defined by the right hand side of
\eqref{5.1a} determines a linear continuous operator on $\SG^1$, i.e., for
arbitrary $p\ge 1$ and $k\ge k_0$, there are $p'\ge p$ and $k'\ge k$ such that
$A_x$ acts continuously from $\SG_{p',k'}^1$
into $\SG_{p,k}^1$. On the other hand, $\di_x$ acts also continuously on
$\SG^1$, and in particular on each $\SG^1_{p',k'}$. Hence, it suffices to prove
the equality \eqref{5.1a} on a total set in $\SG_{p',k'}^1$. As such a set we
take the functions $\wexp{\cdot}{\ttt}$ with $|\ttt|_{p'}<2^{-k'/2}$. Then,
\[ \di_x\wexp{\om}{\ttt}=\ttt(x)\wexp{\om}{\ttt},\]
   and from \eqref{silva} we derive that
\begin{align}
\int_0^\infty\wexp{\om+s\delta_x}{\ttt}\,e^{-s}\,ds&=
\int_0^\infty \exp\big[
\la\om+s\delta_x,\frac\ttt{\ttt+1}\ra-\la\log(1+\ttt)\ra
\Big]
e^{-s}\,ds\notag\\
&=\exp\Big[
\la\om,\frac{\ttt}{\ttt+1}\ra-\la\log(1+\ttt)\ra
\Big]
\int_0^\infty\exp\Big[
\frac{s\ttt(x)}{\ttt(x)+1}-s
\Big]\,ds\notag\\
&=\wexp{\om}{\ttt}\int_0^\infty\exp\Big[
\frac{-s}{\ttt(x)+1}
\Big]\,ds\notag\\
&=\wexp\om\ttt(\ttt(x)+1),
\label{syr}\end{align}
which proves \eqref{5.1a}. The formula \eqref{5.2a} can be proved absolutely
analogously.
\qquad$\blacksquare$
\vspace{2mm}

\noindent{\it Remark\/} 9. Let $\muCP$ be a measure of compound Poisson (CP)
white noise on $S'$ with L\'evy measure $\sigma\nu$, i.e., the Laplace transform
of $\muCP$ is given by
\[\int_{S'}\exp[\la\om,\varphi\ra]\,d\muCP(\om)=\exp\Big[\int_{\R^{d+1}}(e^{s\varphi(x)}-1)\,d\sigma(x)\,d\nu(s)\Big].\]
Then, one can study CP analysis by using an isomorphism between the usual Fock
space over $L^2(\R^{d+1},d\sigma\,d\nu)$ and the $L^2$-space of CP white noise
$(L^2_{\mathrm CP})=L^2(S',{\cal B}(S'),d\muCP)$
\cite{LyReShch,Silva,compound}.

The CP-field operators $a(\xi)$, $\xi\in S$, have in the Fock space the
following representation:
\[a(\xi)=a^+(\xi\ot\id)+a^0(\xi\ot\id)+\int_{\R^d}\xi(x)\,d\sigma(x)\int_\R
s\,d\nu(s)\id+a^-(\xi\ot\id),\]
where $\xi\ot\id=(\xi\ot\id)(x,s)=\xi(x)s$ and $a^+$, $a^0$, $a^-$ are the
standard creation, neutral, and annihilation operators, respectively.

(The image of) the annihilation operator $a^-(\xi\ot\id)$ acts in the following
way in the space $(L^2_{\mathrm CP})$:
\begin{equation}\label{an}
\big(
a^-(\xi\ot \id)\vvv
\big)(\om)=\int_{\R^{d+1}}\xi(x)s(\vvv(\om+s\delta_x)-\vvv(\om))\,d\sigma(x)\,d\nu(s).\end{equation}

In case of the Gamma measure, $d\nu(s)=\frac{e^{-s}}{s}\,ds$ (see \eqref{CP}),
and therefore  the formula \eqref{an} takes the form
\[ \big(
a^-(\xi\ot\id)\vvv
\big)(\om)=\int_{\R^{d+1}}\xi(x)\big(\vvv(\om+s\delta_x)-\vvv(\om)\big)e^{-s}\,d\sigma(x)\,ds.\]
Hence, $a^-(\xi\ot\id)$ is just the operator $a_1^-(\xi)$ under consideration in
this paper.

As a consequence of this, we can write down the explicit action of the adjoint
operator of $a_1^-(\xi)$ in $(L^2_{\mathrm G})$, denoted by $a_1^+(\xi)$, which corresponds to the
operator $a^+(\xi\ot\id)$, see \cite{NuVi,Silva}.
Thus, for an arbitrary $\om\in S'$ of the form \eqref{configuration}
\[
a_1^+(\xi)\phi(\om)=\la\om(x),\xi(x)\phi(\om-s_x\delta_x)\ra-\la\xi\ra\phi(\om),
\qquad \phi\in\operatorname{Dom} a^+(\xi).\]

\section{Creation, neutral, and Gamma annihilation \\ operators on the Gamma space}

In this section, we will obtain the explicit formulas for the operators
$a^+(\xi)$, $a^0(\xi)$, and $a^-_2(\xi)$ acting on $\ltwoG$ by using the formula
\eqref{multiplication}, which expresses the operator of coordinate multiplication
via $\di_x$ and $\dig_x$ and Theorem~6.

Let us fix arbitrary $\tilde p\ge1$ and $\tilde k\ge k_0$, then each $\di_x$ can
be extended to a continuous operator on $\SG^1_{\tilde p,\tilde k}$ and $\dig_x$
and $\om(x)\cdot$ to continuous operators from $\SG_{\tilde p,\tilde k}^1$ into
$\SG_{-\tilde p,-\tilde k}^{-1}$. Choose now $p\ge \tilde p$ and $k\in\N_0$ so
that the space ${\cal E}_{p,k}^1(S')$ is topologically embedded into
$\SG^1_{\tilde p,\tilde k}$.
The restrictions  of $\di_x$, $\dig_x$, and
$\om(x)\cdot$ to ${\cal E}_{p,k}^1(S')$ are continuous operators from this space
into $\SG_{\tilde p,\tilde k}^1$ and $\SG_{-\tilde p,-\tilde k}^{-1}$,
respectively. We rewrite the formula \eqref{multiplication} as follows:
\begin{equation}\label{nov}\om(x)\cdot=\dig_x(\di_x+1)^2+(\di_x+1).\end{equation}
Due to the (proof of) Theorem~6, we have for an arbitrary $\phi\in{\cal
  E}_{p,k}^1(S')$:
\[(\di_x+1)\phi(\om)=\int_0^\infty\phi(\om+s\delta_x)e^{-s}\,ds,\]
and hence on the total set in ${\cal E}_{p,k}^1(S')$ consisting of the functions
 $e^{\la\om,\ttt\ra}$ with $|\ttt|_p< 2^{-k/2}$ we have analogously to
 \eqref{syr}:
\[(\di_x+1)e^{\la \om,\ttt\ra}=(1-\ttt(x))^{-1}e^{\la\om,\ttt\ra}.\]
Therefore, by \eqref{nov},
\[ \om(x)\cdot e^{\la\om,\ttt\ra}=(1-\ttt(x))^{-2}\dig_x
e^{\la\om,\ttt\ra}+(1-\ttt(x))^{-1}e^{\la\om,\ttt\ra},\]
which yields
\begin{align*}
\di_x^+ e^{\la\om,\ttt\ra}&=\om(x)\cdot (\ttt(x)-1)^2
e^{\la\om,\ttt\ra}+(\ttt(x)-1)e^{\la\om,\ttt\ra}\\
&=\om(x)\cdot(\nabla_x-1)^2 e^{\la\om,\ttt\ra}+(\nabla_x-1)e^{\la\om,\ttt\ra}.
\end{align*}
Since $\nabla_x$ is a continuous operator on ${\cal E}_{p,k}^1(S')$, this
immediately implies the central lemma of this section:

\begin{lem} We have on $\SG^1$\rom:
\[\dig_x=\om(x)\cdot(\nabla_x-1)^2+(\nabla_x-1),\]
\end{lem}

Now, by using \eqref{4.5} and Lemma 4, we can calculate, at least formally, the
action of $a^+(\xi)$:
\begin{align*}
a^+(\xi)&=\int_{\supp\sigma}d\sigma(x)\,\xi(x)\big(\om(x)\cdot(\nabla_x-1)^2+(\nabla_x-1)\big)\\
&=\la\om(x),\xi(x)(\nabla_x-1)^2\ra+\la\xi(x),\nabla_x-1\ra.\end{align*}

\begin{th} For any $\phi\in\SG^1$\rom,
\[
a^+(\xi)\phi(\om)=\la\om(x),\xi(x)(\nabla_x-1)^2\phi(\om)\ra+(D_\xi-\la\xi\ra)\phi(\om),\qquad
\om\in S',\]
where $x$ denotes the variable in which the dualization is carried out\rom.
\end{th}

\noindent{\it Proof}. Analogously to the above, we will consider $a^+(\xi)$ as a
continuous operator from ${\cal E}_{p,k}^1(S')$ into $\SG_{\tilde p,\tilde
  k}^1$. Take again $\phi(\om)=\expp$ with $|\ttt|_p<2^{-k/2}$. Then, using
\eqref{4.5} and Lemma~4, we have for an arbitrary $\psi\in\SG^1$:
\begin{align*}
&\ll a^+(\xi)\expp ,\psi(\om)\rr\\
&\qquad
=\int_{\supp\sigma}d\sigma(x)\,\xi(x)\ll\big(\om(x)\cdot(\nabla_x-1)^2+(\nabla_x-1)\big)\expp,\psi(\om)\rr\\
&\qquad= \int_{\supp\sigma}d\sigma(x)\,\xi(x)(\ttt(x)-1)^2\ll\om(x)\cdot
\expp,\psi(\om)\rr\\
&\qquad\quad+\int_{\supp\sigma}d\sigma(x)\,\xi(x)(\ttt(x)-1)\ll
\expp,\psi(\om)\rr\\
&\qquad=\ll \la\om,\xi(\ttt-1)^2\ra\expp,\psi(\om)\rr+\la\xi,\ttt-1\ra\ll
\expp,\psi(\om)\rr\\
&\qquad=\ll
\la\om(x),\xi(x)(\ttt(x)-1)^2\expp\ra+\la\xi,\ttt-1\ra\expp,\psi(\om)\rr.
\end{align*}
Therefore,
\begin{equation}\label{6.1}
  a^+(\xi)\expp=\la\om(x),\xi(x)(\ttt(x)-1)^2\expp\ra+\la\xi,\ttt-1\ra\expp,\end{equation}
where the equality in understood as that in $\SG^{-1}$. But the right  hand side
of \eqref{6.1} considered as a function of $\om$ belongs to ${\cal
  E}_{p,k}^1(S')$, which follows from the representation
\[\la\om(x),\xi(x)(\ttt(x)-1)^2\expp\ra=\sum_{n=1}^\infty\la\om^{\ot n},\frac
1{(n-1)!}\,\ttt^{\ot (n-1)}\ho(\xi(\ttt-1)^2)\ra.\]
Therefore, \eqref{6.1} holds for each $\om\in S_{-p}$.

It remains only to note that $D_\xi$ is a continuous operator on ${\cal
  E}_{p,k}^1(S')$ and for any fixed $\om\in S_{-p}$
  $(\nabla_x-1)^2\phi_m(\om)\to(\nabla_x-1)^2\phi(\om)$ in $S_p$ as a function
  of $x$ if $\phi_m\to\phi$ in ${\cal E}_{p,k}^1(S')$, the latter being proved
  in the same way as the formula \eqref{vier}.\qquad$\blacksquare$

We proceed to consider the neutral operator $a^0(\xi)$ on the Gamma space.

\begin{lem}We have on $\SG^1$\rom:
\[\dig_x\di_x=\om(x)\cdot\nabla_x(1-\nabla_x)-\nabla_x,\qquad x\in\supp\sigma.\]
\end{lem}

\noindent {\it Proof}. Using Lemma 4 and its proof, we get
\begin{gather*}
\dig_x\di_x\expp=\big(\om(x)\cdot(\nabla_x-1)^2+(\nabla_x-1)\big)\di_x\expp\\
=\big(\om(x)\cdot(\nabla_x-1)^2+(\nabla_x-1)\big)\frac{\ttt(x)}{1-\ttt(x)}\,\expp\\
=\bigg(\om(x)\cdot\frac{(\ttt(x)-1)^2\ttt(x)}{1-\ttt(x)}+\frac{(\ttt(x)-1)\ttt(x)}{1-\ttt(x)}\bigg)\expp\\
=\big(\om(x)\cdot(1-\ttt(x))\ttt(x)-\ttt(x)\big)\expp\\
=\big(\om(x)\cdot(1-\nabla_x)\nabla_x-\nabla_x\big)\expp.
\end{gather*}
Again, due to the continuity of $\nabla_x$ on ${\cal E}^1_{p,k}(S')$, we obtain
the lemma.\qquad $\blacksquare$

Since formally
\begin{align*}
a^0(\xi)&=\int_{\supp\sigma}d\sigma(x)\,\dig_x\di_x\\
&=\int_{\supp\sigma}d\sigma(x)\,\xi(x)(\om(x)\cdot\nabla_x(1-\nabla_x)-\nabla_x)\\
&=\la\om(x),\xi(x)\nabla_x(1-\nabla_x)\ra+\la\xi(x),\nabla_x\ra,\end{align*}
we come to the following theorem, whose proof is analogous
to that of Theorem~7.

\begin{th}For any $\phi\in\SG^1$\rom,
\[
a^0(\xi)\phi(\om)=\la\om(x),\xi(x)\nabla_x(1-\nabla_x)\phi(\om)\ra-D_\xi\phi(\om),\qquad
\om\in S'.\]\end{th}

Finally, we will shortly consider the Gamma annihilation operator. Analogously to
Lemma~5, we get

\begin{lem}We have on $\SG^1$
\begin{align*}
\dig_x\di_x^2&=\om(x)\cdot \nabla_x^2-\di_x\nabla_x\\
&=\om(x)\cdot\nabla_x^2-\nabla_x\di_x.\end{align*}
\end{lem}

\begin{th} For any $\phi\in\SG^1$
\begin{align}
a_2^-(\xi)\phi(\om)&=\la\om(x),\xi(x)\nabla_x^2\phi(\om)\ra+D_\xi\phi(\om)\notag\\
&\quad-\int_{\supp\sigma}\int_0^\infty \xi(x)\phi(\om+s\delta_x)e^{-s}\,ds\,d\sigma(x)-\la\xi\ra\phi(\om)
\label{6.2}\\
&=\la\om(x),\xi(x)\nabla_x^2\phi(\om)\ra+D_\xi\phi(\om)\notag\\
&\quad-\int_{\supp\sigma}\int_0^\infty\xi(x)\nabla_x\phi(\om+s\delta_x)e^{-s}\,ds\,d\sigma(x).\label{6.3}
\end{align}\end{th}

\noindent{\it Proof}. The formula \eqref{6.2} follows directly from the equality
\[ a_2^-(\xi)=\la\om,\xi\ra\cdot-a^+(\xi)-2a^0(\xi)-\la\xi\ra\id-a_1^-(\xi)\]
and from Theorems~6--8. The equivalent formula \eqref{6.3} is obtained from
Lemma~6 and Theorem~6 just as above. The only new point here is to prove that
\begin{equation}
\int_{\supp\sigma}d\sigma(x)\,\xi(x)\di_x\nabla_x\phi_m(\om)\to\int_{\supp\sigma}d\sigma(x)\,\xi(x)
\di_x\nabla_x\phi(\om),\qquad \om\in S',\label{6.4}\end{equation}
if $\phi_m\to\phi$ in $\SG^1$ as $m\to\infty$. Bur for each $x\in\supp\sigma$
the norm of $\nabla_x$ on the space ${\cal E}^1_{p,k}(S')$
does not exceed $\|\delta\|_{p,\infty}2^{-k/2}<1$ under an appropriate choice of
$k$, so that each $\nabla_x\di_x=\sum_{n=2}^\infty\nabla_x^n$ (see \eqref{5.1})
is a continuous operator on ${\cal E}_{p,k}^1(S')$ whose norm is bounded by a
constant uniformly in $x$. Therefore, since $\xi\in S\subset
L^1(\supp\sigma,d\sigma)$, $\int_{\supp\sigma}d\sigma(x)\,\xi(x)\nabla_x\di_x$
is a continuous operator on ${\cal E}_{p,k}^1$, which implies \eqref{6.4}.\qquad
$\blacksquare$

\begin{center}\bf ACKNOWLEDGMENTS\end{center}

We would like to thank J. L. Silva for useful
discussions.
The authors were partially supported by the SFB 256, Bonn University.
Yu.K. acknowledges partial financial  support through the INTAS-Project
Nr.~97-0378.

\end{document}